\documentclass{cimart}

\allowdisplaybreaks[4]

\title{On quaternionic analysis and a certain generalized fractal \texorpdfstring{$\psi$}{}-Fueter operator}

\authors{Juan Ram\'irez-Belman, Jos\'e Gonz\'alez-Cervantes and Juan Bory-Reyes}

\authorinfo[J. Ram\'irez-Belman]{Instituto Polit\'ecnico Nacional, M\'exico}{adrianrmzb@gmail.com}

\authorinfo[J. Gonz\'alez-Cervantes]{Instituto Polit\'ecnico Nacional, M\'exico}{jogc200678@gmail.com}

\authorinfo[J. Bory-Reyes]{Instituto Polit\'ecnico Nacional, M\'exico}{juanboryreyes@yahoo.com}

\abstract{%
    This paper introduces a fractal $\psi$-Fueter operator in the quaternionic context inspired by the concepts of proportional derivative and Hausdorff derivative ({fractal derivative}) of a function with respect to a fractal measure. Moreover, we establish the corresponding Stokes and Borel--Pompeiu formulas associated to this generalized fractal $\psi$-Fueter operator with respect to a vector-valued function induced by  truncated exponential functions.
    }


\VOLUME{34}
\YEAR{2026}
\ISSUE{2}
\NUMBER{15}
\DOI{https://doi.org/10.46298/cm.18472}
\editinfo{June 9, 2026}{October 23, 2026}
{Tiago Macedo, Jaqueline Mesquita, Mariel Saez and Rafael Potrie}

\acknowledgments{The work was partially supported by Instituto Polit\'ecnico Nacional (grant numbers IND-2026-0101, SIP20260491) and SECIHTI (grant number 1077475). The authors declare that they have no competing interests regarding the publication of this paper. Finally, the authors thank the anonymous referees for their careful reading of the manuscript and the numerous corrections and suggestions that substantially improved the paper.}

\begin{document}

\section*{Introduction}
\addcontentsline{toc}{section}{Introduction}

The fractal derivative, or Hausdorff derivative, is a relatively new concept of differentiation that extends Leibniz's derivative for discontinuous fractal media. In the literature, there are various definitions of this new concept. For instance, in 2006 Chen introduced the concept of the Hausdorff derivative of a function with respect to a fractal measure ${t^\eta}$, where ${\eta}$ is the order of the fractal derivative, as detailed in \cite{C}. A treatment of a more general case goes back to the work of Jeffery in 1958; see \cite{Je}.

Fractal calculus is extremely effective in branches such as fluid mechanics where hierarchical or porous media, turbulence, or aquifers present fractal properties, which do not necessarily follow a Euclidean geometry.

 {Fractional calculus deals with the generalization of the concepts of differentiation and integration of non-integer orders.} This generalization is not merely a purely mathematical curiosity, but it has demonstrated its application in various disciplines such as physics, biology, engineering, and economics. For a historical review of the theory, we refer the reader to \cite{MR, R, MKM}.

 {Unlike fractional calculus, fractal calculus maintains the chain rule} in a very direct way, which relates the fractal derivative to the classical derivative.

Fractal derivatives form a new class of derivatives with many applications in real-world problems:  hydrodynamics \cite{BE}, fractal viscoelastic models \cite{CCX}, models for suspended-sediment transport in steady flows \cite{NSLZX}. The Hausdorff fractal derivative model has been used to understand complex dynamics in natural systems, including anomalous diffusion; see \cite{LWMag,SLZC}. 

Fractals are complex geometric patterns that repeat at different scales, while fractional derivatives are a generalization of ordinary derivatives that allow for non-integer orders. 

A considerable literature has grown up around   fractal  derivatives.  For references connected with the subject being considered in this  work, we refer the reader to     \cite{AS,C,CSZ,DMP,KH,H,SAH,SAH2,SHL,SSMH}   and the references quoted there.  

Quaternionic analysis (the most natural and close generalization of complex analysis) concerns the connection between analysis (even topology / geometry) in ${\mathbb R^4}$ and the algebraic structure of quaternions ${\mathbb H}$. At the heart of this function theory lies the notion  of ${\psi}$-hyperholomorphic functions defined on domains in ${\mathbb R^4}$ with values in ${\mathbb H}$, i.e., null solutions of the so-called ${\psi}$-Fueter operator (to be defined later) in which the standard basis of ${\mathbb H}$ is replaced by a structural set  ${\psi=\{1, \psi_1, \psi_2, \psi_3\}}$.

{In recent years}, there has been increasing interest in finding a framework for a  fractal  counterpart of quaternionic analysis; see \cite{BPP, GB1, GBS, GBNS} and the references therein. For example, a theory of the quaternionic right-module of fractional hyperholomorphic functions induced by a fractional differential Fueter-type operator, in the Riemann--Liouville and Caputo senses,  was introduced  in \cite{GB1}, and these results were extended in \cite{GB2} in which the concept 
 ``fractional hyperholomorphic functions with respect to  a vector-valued function'' was introduced. 

In addition, it is well known that quaternionic slice regular function theory is  another  natural extension of the holomorphic function theory. Thus, a {theory of the fractional slice regular functions of a quaternionic variable induced by  a fractional Cauchy--Riemann operator,} on each slice, in the Riemann--Liouville and Caputo senses  was presented  in \cite{GBS}.

On the other hand, among the different types of generalized {derivatives} is the quantum derivative. In this setting, the quantum Fueter operator and the   quantum Cauchy--Riemann operator were defined to introduce the $(q,q')$-model of quaternionic analysis and the theory of slice regular functions via post-quantum calculus theory, respectively; see \cite{GBS2, GBNS}.

Until now the theory of ${\mathbb{H}}$-valued   {fractal  hyperholomorphic functions  and  the theory of  fractal slice regular  functions } have grown to be able to establish the ${\mathbb H}$-modules of Bergman and Dirichlet {associated} to fractal Fueter and Cauchy--Riemann operators.

This {work} {introduces} a  fractal ${\psi}$-Fueter operator in the quaternionic context inspired by the concepts of proportional   Hausdorff derivative of a function with respect to a fractal measure. The {${\psi}$-fractal  proportional} hyperholomorphic function theory is also established by extending the well-known results of the ${\psi}$-hyperholomorphic function theory  presented in \cite{Sh,SV1,SV2}.

Moreover, we establish the corresponding Stokes and Borel--Pompeiu formulas associated with this generalized {fractal ${\psi}$-Fueter operator.}

It is important to comment that we use an arbitrary structural set ${\psi}$ because in the quaternionic framework literature, it is well known that many particular structural sets give us interesting relationships between some important function theories and the remarkable ${\psi}$-hyperholomorphic function theory. For example, all different ways in which the theory of holomorphic functions in two complex variables is embedded in the ${\psi}$-hyperholomorphic function theory are given by 
${\psi = \{ 1, i, ie^{i\theta}  j, e^{i\theta} j\},}$
 where ${\theta\in (0,2\pi]}$, see \cite{MS}.

Other important relationships involving particular $\psi$-hyperholomorphies concern the Cimmino system, the Div-Rot system, the Riesz system, and the system for Dirac's bispinor in harmonic time.

The {structure} of this {work} is summarized as follows. In Section \ref{1} we give a brief exposition of the generalized  {fractal  derivative considered}. Section \ref{sec:quaternionic} presents some preliminaries on quaternionic analysis. In Section \ref{sec:general} we develop the rudiments of a function theory induced by a quaternionic ${\beta}$-proportional fractal Fueter operator and finally in Section \ref{sec:truncated} we consider a quaternionic ${\beta}$-proportional fractal Fueter operator with truncated exponential functions {as a fractal} measure.

\section{Generalized fractal derivative}\label{1}

\begin{definition}\label{fractal} The {Hattaf} fractal derivative of a function ${f}$, defined on an interval ${I}$,  with respect to a fractal measure ${\nu(\eta, t)}$ is given  by
\begin{align*} \frac{d_\nu f(t)}{
dt^\eta} := \lim_{\tau \to t} 
\frac{
f(t) - f(\tau)}
{\nu(\eta, t) - \nu (\eta, \tau)} 
, \quad  \eta > 0. 
\end{align*}
If 
\[
\dfrac{d_\nu f(t)}{dt^\eta},
\]
 exists for all ${t \in   I}$, then  ${f}$ is real fractal differentiable on ${I}$ with order ${\eta}$, see \cite{KH}. 
\end{definition}
We record some well-known cases. If ${\nu(h, t) = t}$ for all ${t\in I}$, then 
\[
\dfrac{d_\nu }{dt^\eta} = \dfrac{d }{dt},
\]
 is the derivative operator. In addition, if ${\nu(h, t) = h(t)}$ for all ${t\in I}$ where ${h'(t) > 0}$ for all ${t\in I}$, then 

\[
\dfrac{d_\nu f(t)}{dt^\eta}  =  \dfrac{f'(t)}{h'(t)},
\]
 for all ${f\in C^1(I)}$. {Remember that ${h'(t)\neq 0}$, since our hypothesis is that ${h'(t) > 0}$ for all ${t\in I}$.}

On the other hand, if ${\nu (\eta, t) = t^\eta}$ for all ${t\in I}$, then   
\[
\dfrac{d_\nu f(t)}{dt^\eta},
\]
 reduces to the Hausdorff derivative. Another useful fractal measure is ${\nu (\eta, t)= e^{t^\alpha}}$, for all ${t\in I}$ and ${\alpha \in (0,1]}$.

We consider a well-known extension of the previous fractal derivative.

\begin{definition}\label{defbfractal} Given {${\beta\in (0,1]}$}, we present the ${\beta}$-fractal derivative of a function ${f}$, defined on an interval ${I}$,  with respect to a fractal measure ${\nu(\eta, t)}$:  
\begin{align*} \frac{d^{\beta}_\nu f(t)}{
dt^\eta} := \lim_{\tau \to t} 
\frac{
(f(t))^{\beta} - (f(\tau))^{\beta}}
{\nu(\eta, t) - \nu (\eta, \tau)} 
, \quad  \eta > 0. 
\end{align*}
\end{definition}

In order to make our description of the concept of  {fractal  derivatives to be used precise,} we introduce the notion of fractional proportional derivative, following \cite{JAR}.
\begin{definition}\label{defchi}
For ${\sigma\in[0,1]}$, let ${ \chi_0 ,  \chi_1: [0, 1]\times I }$ be continuous functions such that  
\[\lim_{\sigma\to 0^+} \chi_1(\sigma,t) = 1,  \lim_{\sigma\to 0^+} \chi_0(\sigma,t) = 0, \lim_{\sigma\to 1^-} \chi_1(\sigma,t) = 0, \lim_{\sigma\to 1^-} \chi_0(\sigma,t) = 1,\]
and ${\chi_{1}(\sigma,t)\neq 0}$, ${\sigma\in[0,1)}$, ${\chi_{0}(\sigma,t)\neq 0}$, ${\sigma\in(0,1]}$. Then the proportional derivative of ${f\in C^1(I)}$ of order ${\sigma}$ is given by
\begin{align*} 
D^{\sigma}f (t) =  \chi_1(\sigma,t)f (t) +  \chi_0(\sigma,t)f'(t) , 
\end{align*}
{for all ${t\in I}$.}
\end{definition}
  
Combining Definitions \ref{defbfractal} and \ref{defchi} yields the following definition.
  
\begin{definition}\label{defpropbfractal} 
Let ${\beta\in (0,1]}$. The proportional ${\beta}$-fractal derivative of ${f:I\to\mathbb R}$ with respect to ${\nu(\eta, t)}$ and ${\sigma}$ is defined to be 
\begin{align*}
\frac{d^{\sigma,\beta}_\nu f(t)}{
dt^\eta} :=  \chi_1(\sigma,t)f (t) +  \chi_0(\sigma,t)
\frac{d^{\beta}_\nu f(t)}{dt^\eta}, 
\end{align*}
if it exists for all ${t\in I}$.
\end{definition}

\begin{remark}\label{rem1}
Given ${\alpha \in (0,1]}$ and {${k\in \mathbb N:=\{1,2,3,\dots\}}$}, we will consider the ${k}$-truncated exponential function defined as follows
\[\displaystyle e(t^{\alpha})_k := \sum_{i=0}^k \frac{(t^\alpha)^i}{i!},\] 
for all ${t>0}$. For ${k=1}$ we have ${e(t^{\alpha})_1= 1+ t^{\alpha}}$  and  for ${k=\infty}$ we {denote}  ${e(t^{\alpha})_\infty = e^{t^{\alpha}}}$. 
\end{remark}

\begin{remark}\label{rem2} An important particular case is obtained when the proportional coefficients and fractal measure in Definition \ref{defpropbfractal} are given by 
\[\chi_1(\sigma,t) = 1-\sigma , \quad 
\chi_0(\sigma,t)=\sigma, \quad\nu(k, t)= e(t^{\alpha})_k,
\] for all ${\sigma\in [0,1]}$ and ${t\in I}$.
This choice gives the following generalized fractal derivatives. For ${f\in C^1(I)}$ we have 
\begin{align*}  
\frac{ {d}^{\sigma,\beta}f(t)}{ d t _{\alpha, k} }: = (1-\sigma)f (t) +   \sigma \frac{ (f ^{\beta}) '(t) }{   e ( t^{\alpha})_{k }'},
\end{align*}
for all ${t\in  I}$. The particular cases ${k=1,\infty}$ reduce to
\begin{align*}  
\frac{ {d}^{\sigma,\beta}f(t)}{ d t _{\alpha, 1} } &= (1-\sigma)f(t) + \sigma \frac{ (f ^{\beta}) '(t)}{\alpha t^{\alpha-1}},\\
 \displaystyle \frac{ {d}^{\sigma,\beta}f(t)}{d t_{\alpha, \infty} } &= (1-\sigma)f (t) + \sigma \frac{(f ^{\beta})'(t)}{\alpha t^{\alpha-1} e^{t^{\alpha}}},
\end{align*}
where the conditions ${\alpha\in (0,1]}$ and ${t>0}$ are necessary.

Addressing the issue ${\sigma=\alpha}$ requires  that the case ${\alpha=0}$ should be omitted.   
\begin{align*}  
\frac{{d}^{\alpha,\beta}f(t)}{d t _{\alpha, 1}} &= (1-\alpha)f (t) + \frac{ (f ^{\beta}) '(t)}{t^{\alpha-1}},\\
\displaystyle \frac{{d}^{\alpha,\beta}f(t)}{d t _{\alpha, \infty}} &= (1-\alpha)f (t) + \frac{ (f ^{\beta}) '(t)}{t^{\alpha-1} e^{t^{\alpha}}}.
\end{align*}

The ${k}$-truncated exponential function as fractal measure provides the generalized  {fractal  derivative} in the following general form:
\begin{align*}  
\frac{{d}^{\alpha,\beta}f(t)}{d t _{\alpha, k} } := (1-\alpha)f (t) + \frac{ (f ^{\beta}) '(t)}{t^{\alpha-1} e(t^{\alpha} )_{k-1}}.
\end{align*}
\end{remark}

\section{Preliminaries on quaternionic analysis}\label{sec:quaternionic}
Nowadays, quaternionic analysis is regarded as a broadly accepted branch of classical analysis offering a successful generalization of complex holomorphic function theory. The most renowned works are \cite{Fu1, Fu2, S} and the books \cite{GS1, GS2, KS, K}. It relies heavily on results on functions defined on domains in $\mathbb R^4$ with values in the skew field of real quaternions $\mathbb H$ associated with a generalized Cauchy--Riemann operator (the so-called $\psi$-Fueter operator) by using a general orthonormal basis in $\mathbb R^4$ called a structural set.

We begin by recalling some background and fixing notation that will be used throughout the entire document. For more details, we refer the interested reader to \cite{SV1, SV2, S}.

A real quaternion is an element of the form 
\[
x=x_0  + x_{1} {e_1} +x_{2} e_2 + x_{3} e_3,
\]
 where ${x_0, x_1, x_2, x_3\in\mathbb R}$ and  the imaginary units ${e_1}$, ${e_2}$, ${e_3}$ satisfy 

\[e_1^2=e_2^2=e_3^2=-1,  e_1e_2=-e_2e_1=e_3, e_2e_3=-e_3e_2=e_1, e_3e_1=-e_1e_3=e_2.\]

Quaternions form a skew field denoted by  ${\mathbb H}$. The set ${\{1,e_1,e_2,e_3\}}$ is the standard basis of ${\mathbb H}$.

The vector part of ${x\in \mathbb H}$ is, by definition, 
\[
{\bf{x}}:= x_{1} {e_1} +x_{2} e_2 + x_{3} e_3,
\]
 while its real part is ${x_0:=x_0}$.
The quaternionic conjugation of ${x}$, denoted by ${\bar x}$, is defined by ${\bar x:=x_0-{\bf x}}$ and the norm of ${x\in\mathbb H}$ is given by
\[\|x\|:=\sqrt{x_0^2 +x_1^2+x_2^2+x_3^2}= \sqrt{x\bar x} = \sqrt{\bar x  x}.\]

The quaternionic scalar product  of ${x, y\in\mathbb H}$ is given by 
\[\langle x, y\rangle:=\frac{1}{2}(\bar x y + \bar y x) = \frac{1}{2}(x \bar y + y  \bar x).\]

A 4-tuple  of quaternions  ${\psi=(\psi_0, \psi_1,\psi_2,\psi_3)}$ is called a structural set if ${\langle \psi_k, \psi_s\rangle =\delta_{k,s}}$, 
for ${k, s=0,1,2,3}$. Any quaternion ${x}$  can be rewritten as ${x_{\psi} := \sum_{k=0}^3 x_k\psi_k}$, where ${x_k\in\mathbb R}$ for all ${k}$. Basic properties of structural sets can be found in \cite{SV1,SV2,MS,Sh}.

Given ${q, x\in \mathbb H}$, we follow the notation used in {\cite{SV1}} to write
\[\langle q, x \rangle_{\psi}=\sum_{k=0}^3 q_k x_k,\] 
where ${q_k, x_k\in \mathbb R}$ for all ${k}$. 

Given a domain ${\Omega\subset\mathbb H\cong \mathbb R^4}$  and a function ${f:\Omega\to \mathbb H}$, we write ${f}$ as ${f=\sum_{k=0}^3 f_k \psi_k}$, where ${f_k, k= 0,1,2,3,}$ are ${\mathbb R}$-valued functions. Properties of ${f}$, including continuity, differentiability, and integrability, are defined in terms of the corresponding properties of its components ${f_k}$. For example, ${C^{1}(\Omega, \mathbb H)}$ denotes the set of continuously differentiable ${\mathbb H}$-valued functions defined in ${\Omega}$.  

The left- and right-${\psi}$-Fueter operators are given by   

\begin{align*}
{}^{{\psi}}\mathcal D[f] :=& \sum_{k=0}^3 \psi_k \partial_k f,\\
\intertext{and}
{}^{{\psi}}\mathcal  D_r[f] :=&  \sum_{k=0}^3 \partial_k f \psi_k
,
\end{align*}
 for all ${f \in C^1(\Omega,\mathbb H)}$, respectively, where 
\[
\partial_k f =\displaystyle \frac{\partial f}{\partial x_k},
\]
 for all ${k}$. 

When ${\psi=\{1,e_1,e_2,e_3\}}$ the ${\psi}$-Fueter operators reduce to standard Fueter operators in classical quaternionic analysis. 

Let ${\partial \Omega}$ be a 3-dimensional compact smooth surface. Recall the Borel--Pompeiu formula and the differential and integral versions of Stokes' formula:  
\begin{align}\label{BorelHyp}  &  \int_{\partial \Omega}(K_{\psi}(\tau-x)\sigma_{\tau}^{\psi} f(\tau)  +  g(\tau)   \sigma_{\tau}^{\psi} K_{\psi}(\tau-x) ) \nonumber  \\ 
&  - 
\int_{\Omega} (K_{\psi} (y-x) {}^{\psi}\mathcal D [f] (y) + {}^{{\psi}}\mathcal  D_r [g] (y) K_{\psi} (y-x)
     )dy   \nonumber \\
	=  &  \left\{ \begin{array}{ll}  f(x) + g(x) , &  x\in \Omega,  \\ 0 , &  x\in \mathbb H\setminus\overline{\Omega},                     
\end{array} \right. 
\end{align}  
for all ${f,g\in C^{1}(\overline{\Omega}, \mathbb H)}$. 

\begin{align}\label{StokesHyp} d(g\sigma^{{\psi} }_x f) = & \left(g \, {}^{{\psi}}\mathcal  D[f]+ \  {}^{{\psi}}\mathcal D_r[g] \,f\right)dx, \nonumber\\
\int_{\partial \Omega} g\sigma^\psi_x f =  &   \int_{\Omega } \left( g\, {}^\psi \mathcal  D[f] + {}^{{\psi}}\mathcal  D_r[g] \,f\right)dx,
\end{align}
for all ${f,g \in C^1(\overline{\Omega}, \mathbb H)}$. Here ${d}$ represents the exterior differentiation operator, ${dx}$ is  the differential form of the 4-dimensional volume in ${\mathbb R^4}$ and  
\[\sigma^{{\psi} }_{x}:= sgn\psi  \sum_{k=0}^3 (-1)^{k+1} \psi_k d\hat{x}_k,\]
is the quaternionic differential form of the 3-dimensional volume in ${\mathbb R^4}$ according  to ${\psi}$, where 
\[
d\hat{x}_k  = dx_0 \wedge dx_1\wedge dx_2  \wedge  dx_3,
\]
 omitting factor ${dx_k}$.  {In addition,} ${sgn\psi}$ is ${1}$, or ${-1}$,  if  ${\psi}$ and  
\[
\psi_{std}:=\{1, {\bf i}, {\bf j}, {\bf k}\},
\]
 have the same orientation, or not, respectively.  Note that ${|\sigma^{{\psi} }_{x}| = dS_3}$  is the differential form of the 3-dimensional volume in ${\mathbb R^4}$ and write ${\sigma_x=\sigma^{{\psi_{std}} }_{x}}$. Let us recall that the ${\psi}$-Cauchy kernel is given by 
 \[ K_{\psi}(\tau- x)=\frac{1}{2\pi^2} \frac{ \overline{\tau_{\psi} - x_{\psi}}}{|\tau_{\psi} - x_{\psi}|^4}.\]

\section{Function theory generated by a  proportional fractal Fueter operator}\label{sec:general}
Let us extend Definition \ref{defpropbfractal} to a quaternionic  differential operator associated with an arbitrary structural set ${\psi}$. Let ${\Omega\subset \mathbb H}$ be a domain  from now on.

\begin{definition}\label{Hdefpropbfractal}
Fix
\[
\beta=(\beta_0,\beta_1,\beta_2,\beta_3) \in (0,1]^4,
\]
and ${\nu= (\nu_0, \nu_1,\nu_2,\nu_3)}$ where ${\nu_k(\eta_k, x_k)}$ is a fractal measure for ${k=0,1,2,3}$ according to Definition \ref{fractal}. 
Denote 

\[
\chi_1= (\chi_{0,1}, \chi_{1,1} ,\chi_{2,1},\chi_{3,1} )
,\]

\[
\chi_0= (\chi_{0,0}, \chi_{1,0} ,\chi_{2,0},\chi_{3,0} ),
\]
 and 
\[
\sigma=(\sigma_0, \sigma_1,\sigma_2,\sigma_3) \in [0,1]^4
,\]
where ${\chi_{k,1}(\sigma_k, x_k)}$ and ${\chi_{k,0}(\sigma_k, x_k)}$ are given by Definition \ref{defchi} on coordinate ${x_k}$, ${k=0,1,2,3}$.

Denote 
\[
\mathbb H^+:=\{x= x_0+ x_1e_1+ x_2e_2 + x_3e_3 \in \mathbb H \ \mid \  x_0,x_1,x_2,x_3 >0\}
.\]
 In this subsection,  we consider functions ${f :\Omega\to \mathbb H}$ satisfying ${f(\Omega) \subset \mathbb H^{+}}$ and assume that 
\[
\dfrac{\partial^{\beta_n}_{\nu_n} f(x)}{\partial (x_n)^{\eta_n}},
\]
 exists for all ${x\in \Omega}$ and all ${n=0,1,2,3}$. Then, the (left) quaternionic ${\psi}$-proportional ${\beta}$-fractal Fueter-type operator of ${f}$  with respect to ${\nu}$ and ${\sigma}$, is given by    
\begin{align*}
{}^{\psi}\mathcal D^{\sigma,\beta}_{\nu} [f] (x)
:= &
\sum_{n=0}^3 \psi_n \frac{\partial^{\sigma_n,\beta_n}_{\nu_n} f(x)}{
\partial (x_n)^{\eta_n} }  \\
 = &  \sum_{n=0}^3 \psi_n \left( \chi_{n,1}(\sigma_n,x_n )f (x) +  \chi_{n,0}(\sigma_n, x_n)
\frac{\partial^{\beta_n}_{\nu_n} f(x)}{
\partial (x_n)^{\eta_n}} \right)  \\
= &  \sum_{n, m = 0}^3 \psi_n\psi_m \left( \chi_{n,1}(\sigma_n,x_n )f_m (x) +  \chi_{n,0}(\sigma_n, x_n)
\frac{\partial^{\beta_n}_{\nu_n} f_m(x)}{\partial (x_n)^{\eta_n}} \right).  
\end{align*}
The set  of (left) ${\sigma}$-proportional, ${\beta}$-fractal ${\psi}$-hyperholomorphic  functions with respect to ${\nu}$ defined on ${\Omega}$ is denoted by 
\[
{}^{\psi}\mathcal M^{\sigma,\beta}_{\nu}(\Omega, \mathbb H),
\]
 and consists of ${f\in C^1(\Omega, \mathbb H)}$  such that
\[{}^{\psi}\mathcal D^{\sigma,\beta}_{\nu} [f] (x)=0,\] {for all ${x\in\Omega}$.}
\end{definition}

\begin{definition}\label{Right-Hdefpropbfractal}
Fix
\[
\delta=(\delta_0,\delta_1,\delta_2,\delta_3)  \in (0,1]^4,
\]
and ${\mu= (\mu_0, \mu_1,\mu_2,\mu_3)}$ where ${\mu_k(\zeta_k, x_k)}$ is a fractal measure for ${k=0,1,2,3}$ according to Definition \ref{fractal}. 
Denote 
\[
\varkappa_1= (\varkappa_{0,1}, \varkappa_{1,1} ,\varkappa_{2,1},\varkappa_{3,1})
,\]
 
\[
\varkappa_0= (\varkappa_{0,0}, \varkappa_{1,0} ,\varkappa_{2,0},\varkappa_{3,0}),
\]
 and 
\[
\rho=(\rho_0, \rho_1,\rho_2,\rho_3) \in [0,1]^4,
\]
 where ${\varkappa _{k,1}(\rho_k, x_k)}$ and ${\varkappa_{k,0}(\rho_k, x_k)}$ are given by Definition \ref{defchi} for ${k=0,1,2,3}$.
 
Let ${g:\Omega\to \mathbb H}$ be such that 
\[
\dfrac{\partial^{\delta_n}_{\mu_n} g(x)}{\partial (x_n)^{\zeta_n}},
\]
 exists for all ${x\in \Omega}$ and all ${n=0,1,2,3}$. The (right) quaternionic ${\psi}$-proportional ${\delta}$-fractal Fueter-type operator of ${g}$  with respect to ${\mu}$ and ${\rho}$, is given by    
\begin{align*}
{}^{\psi}\mathcal D^{\rho,\delta}_{r,\mu} [g] (x)
:= &
\sum_{n=0}^3\frac{\partial^{\rho_n,\delta_n}_{\mu_n} g(x)}{
\partial (x_n)^{\zeta_n} }  \psi_n \\
= &  \sum_{n, m = 0}^3 \psi_m\psi_n \left(  \varkappa_{n,1}(\rho_n,x_n ) g_m (x) +  \varkappa_{n,0}(\rho_n, x_n)
\frac{\partial^{\delta_n}_{\mu_n} g_m(x)}{
\partial (x_n)^{\zeta_n}} \right) . 
\end{align*}
The set of (right) ${\rho}$-proportional, ${\delta}$-fractal  ${\psi}$-hyperholomorphic  functions with respect to ${\mu}$ defined on ${\Omega}$ is denoted by 
\[
{}^{\psi}\mathcal M^{\rho,\delta}_{r,\mu}(\Omega, \mathbb H),
\]
 and consists of ${g\in C^1(\Omega, \mathbb H)}$  such that
\[{}^{\psi}\mathcal D^{\rho,\delta}_{r,\mu} [g]  (x)=0,\] {for all ${x\in\Omega}$.}
\end{definition}

\begin{remark}
Note that  
\[
f\in {}^{\psi}\mathcal M^{\sigma,\beta}_{\nu}(\Omega, \mathbb H),
\]
 and 
\[
g\in {}^{\psi}\mathcal M^{\rho,\delta}_{r,\mu} (\Omega, \mathbb H),
\]
 if and only if the proportional-fractal Cauchy--Riemann-type equations are satisfied by the real components of ${f}$ and ${g}$. What is more, for some specific functions ${\sigma}$, ${\beta}$ and ${\nu}$ we can see that 
\[
{}^{\psi}\mathcal M^{\sigma,\beta}_{\nu}(\Omega, \mathbb H),
\]
 is the ${\psi}$-hyperholomorphic function theory studied in \cite{SV1,SV2}.
\end{remark}
The following proposition presents identities in which operators ${{}^{\psi} \mathcal D}$ and ${{}^{\psi}\mathcal D^{\sigma,\beta}_{\nu}}$ are related. Similarly, ${{}^{\psi} \mathcal D_r}$ and ${{}^{\psi}\mathcal D^{\rho,\delta}_{r,\mu}}$ are related.

\begin{definition}
Given ${f, g\in C^1(\Omega,\mathbb H)}$ as above, let us assume that   
\begin{align*}
 \lambda^{\beta_n}_{\nu_n}(f_m) (x):= & \int_{0}^{x_n  } \frac{\partial^{\beta_n}_{\nu_n} f_m(x)}{
\partial (t_n)^{\eta_n}}dt, \quad 
 \lambda^{\beta_n}_{\nu_n}(f ) (x):=   \sum_{m=0}^3 \psi_m\int_{0}^{x_n  } \frac{\partial^{\beta_n}_{\nu_n} f_m(x)}{
\partial (t_n)^{\eta_n}} dt  
,\\
 \lambda^{\delta_n}_{\mu_n}(g_m) (x)=  &  \int_{0}^{x_n  } \frac{\partial^{\delta_n}_{\mu_n} g_m(x)}{
\partial (t_n)^{\zeta_n}}  dt   ,\quad 
 \lambda^{\delta_n}_{\mu_n}(g ) (x)=   \sum_{m=0}^3 \psi_m\int_{0}^{x_n  } \frac{\partial^{\delta_n}_{\mu_n} g_m(x)}{
\partial (t_n)^{\zeta_n}}dt,
\end{align*}
for all ${x\in\Omega}$, exist for ${n,m=0,1,2,3}$. Suppose that  ${\chi_{n,0}}$ and  ${\varkappa_{n,0}}$  are functions of class ${C^1}$ such that  ${\chi_{n,0}(\sigma_n, x_n )\neq 0}$ and ${\lambda^{\beta_n}_{\nu_n}(f_m) (x)\neq 0}$, 
${\varkappa_{n,0}(\rho_n, x_n )\neq 0}$  and ${\lambda^{\delta_n}_{\mu_n} (g_m) (x)\neq 0}$ for all ${x=(x_0,x_1,x_2,x_3) \in \Omega}$ and all  ${m,n=0,1,2,3}$. Denote 
\begin{align*}
\mathcal E^{\sigma,\beta}_{\nu}  [f](x) :=&\sum_{{\substack{\ell, n = 0 \\ \ell \ne n} } }^3 \psi_n     \frac{\partial }{\partial x_n} \left[ ( \chi_{k,0}(\sigma_k, x_k )   )     \lambda^{\beta_k}_{\nu_k}(f ) (x)  \right] ,\\ 
\mathfrak L^{\sigma,\beta}_{\nu}  (f)  (x):= &  \sum_{k=0 }^3  ( \chi_{k,0}(\sigma_k, x_k )   )  \lambda^{\beta_k}_{\nu_k}(f ) (x) , \\
  L_{n,m}[f](x) := & \frac{\partial }{\partial x_n} \left(   \frac{  \chi_{n,0}(\sigma_n,x_n )   } {  e^{h_{n,m}(x) }    }
\right)  e^{h_{n,m}(x) },\\
h_{n,m}(x) := & \int_{0}^{x_n} \frac{\chi_{n,1}(\sigma_n,t_n)}{\chi_{n,0}(\sigma_n,t_n)} \frac{f_m}{\lambda^{\beta_n}_{\nu_n}(f_m)}dt,
\end{align*}
and 
\begin{align*}
\mathcal E^{\rho,\delta}_{r,\mu}  [g](x):=& \sum_{{\begin{array}{c} n=0=k \\  n\neq k \end{array} } }^3     \frac{\partial }{\partial x_n} \left[ ( \varkappa_{k,0}(\rho_k, x_k )) \lambda^{\delta_k}_{\mu_k}(g ) (x)  \right] \psi_n ,  \\
\mathfrak L ^{\rho,\delta}_{\mu} (g)  (x) :=& \sum_{k=0 }^3  ( \varkappa_{k,0}(\rho_k, x_k))  \lambda^{\delta_k}_{\mu_k}(g ) (x), \\
  T_{n,m}[g](x) := & \frac{\partial }{\partial x_n} (\frac{  \varkappa_{n,0}(\rho_n,x_n)} {e^{l_{n,m}(x)}}
)  e^{l_{n,m}(x) },\\
 l_{n,m}(x) := & \int_{0}^{x_n} \frac{ \varkappa_{n,1}(\rho_n,t_n )}{ \varkappa_{n,0}(\rho_n,t_n )}\frac{g_m}{\lambda^{\delta_n}_{\mu_n}(g_m)}dt,
\end{align*}
for all ${n,m\in \{0,1,2,3\}}$.
\end{definition}

\begin{proposition}\label{Fueter-Hdefpropbfractal}
If ${f, g\in C^1(\Omega,\mathbb H)}$ satisfy  the previous  definition, then
\begin{align*}
{}^{\psi} \mathcal D \circ   \mathfrak L ^{\sigma,\beta}_{\nu} [f] (x) =&  {}^{\psi}\mathcal D^{\sigma,\beta}_{\nu} [f] (x) +
\mathcal E^{\sigma,\beta}_{\nu}  [f](x) +
\sum_{{ n, m = 0} }^3 \psi_n\psi_m  L_{n,m}[f](x) \lambda^{\beta_n}_{\nu_n}(f_m) (x),\\
{}^{\psi} \mathcal D_r \circ \mathfrak L ^{\rho,\delta}_{\mu}[g](x) =& {}^{\psi}\mathcal D^{\rho,\delta}_{r,\mu} [g] (x) +
\mathcal E^{\rho,\delta}_{r,\mu}  [g](x) + \sum_{{ n, m = 0} }^3 \psi_m\psi_n  T_{n,m}[g](x) \lambda^{\delta_n}_{\mu_n}(g_m) (x),
\end{align*}
for all ${x\in \Omega}$.
\end{proposition}

\begin{proof} The proofs of the two identities are similar, so we present only the first.

To simplify notation consider ${\lambda_n = \lambda^{\beta_n}_{\nu_n}}$ for all ${n=0,1,2,3}$.  {From direct computations that  involve properties of the exponential and a version of the fundamental theorem of calculus in one real variable,}  we have that 
\begin{align*}
\frac{\partial }{\partial x_n} (e^{h_{n,m}(x)}\lambda_{n}(f_m) (x) ) = { }&{ } e^{h_{n,m}(x) }  \left[ \frac{ \chi_{n,1}(\sigma_n, x_n )}{ \chi_{n,0}(\sigma_n,  x_n)} \frac{f_m}{\lambda_{n}(f_m)}\lambda_{n}(f_m) (x) + \frac{\partial^{\beta_n}_{\nu_n} f_m(x)}{
\partial (x_n)^{\eta_n}} \right] \\
= { }&{ } \frac{ e^{h_{n,m}(x)}}{\chi_{n,0}(\sigma_n,x_n )}\left[\chi_{n,1}(\sigma_n,t_n)f_m + \chi_{n,0}(\sigma_n, x_n)\frac{\partial^{\beta_n}_{\nu_n} f_m(x)}{\partial (x_n)^{\eta_n}} \right] \\
= { }&{ } \frac{ e^{h_{n,m}(x)} } {\chi_{n,0}(\sigma_n,x_n )} \frac{\partial^{\sigma_n,\beta_n}_{\nu_n} f_m(x)}{\partial (x_n)^{\eta_n}}, \\
\frac{\partial^{\sigma_n,\beta_n}_{\nu_n} f_m(x)}{\partial (x_n)^{\eta_n} }  = { }&{ }  \frac{  \chi_{n,0}(\sigma_n,x_n)}{e^{h_{n,m}(x)}} 
\frac{\partial }{\partial x_n} \left(e^{h_{n,m}(x) }  \lambda_{n}(f_m) (x) \right),
\end{align*}
for  ${n=0,1,2,3}$.
Then, 
\begin{align*}
^{\psi}\mathcal D^{\sigma,\beta}_{\nu} & [f] (x) \\
 = { }&{ }   \sum_{n, m =0}^3 \psi_n\psi_m   \frac{  \chi_{n,0}(\sigma_n,x_n )   } {  e^{h_{n,m}(x) }    } 
   \frac{\partial }{\partial x_n} \left(e^{h_{n,m}(x) }  \lambda_{n}(f_m) (x) \right) \\
 = &  \sum_{n, m = 0}^3 \psi_n\psi_m  \frac{\partial }{\partial x_n} (\chi_{n,0}(\sigma_n,x_n) \lambda_{n}(f_m) (x) ))  - 
 \sum_{n, m = 0}^3 \psi_n\psi_m  L_{n,m}[f](x) \lambda_{n}(f_m) (x) \\
 = &  \sum_{n=0 }^3 \psi_n  \frac{\partial }{\partial x_n} (\chi_{n,0}(\sigma_n,x_n) \sum_{m=0}^3  \psi_m \lambda_{n}(f_m) (x) ) ) - 
 \sum_{n, m = 0}^3 \psi_n\psi_m  L_{n,m}[f](x) \lambda_{n}(f_m) (x) \\
 = &  \sum_{n=0 }^3 \psi_n \frac{\partial }{\partial x_n} \left[\chi_{n,0}(\sigma_n,x_n ) \lambda_{n}(f) (x) \right]    - 
 \sum_{n, m = 0}^3 \psi_n\psi_m L_{n,m}[f](x) \lambda_{n}(f_m) (x) \\
 = &  \sum_{n=0 }^3 \psi_n \frac{\partial }{\partial x_n} \left[\varsigma(x) \lambda_{n}(f ) (x) \right]    -
 \sum_{n=0 }^3 \psi_n     \frac{\partial }{\partial x_n} \left[\kappa_n \lambda_{n}(f)(x) \right]  
\\ & - \sum_{n, m = 0}^3 \psi_n\psi_m  L_{n,m}[f](x) \lambda_{n}(f_m) (x), 
\end{align*}
for all ${x\in \Omega}$, where
\[
\varsigma(x)=\sum_{\ell=0}^3 \chi_{\ell,0}(\sigma_\ell , x_\ell),
\]
 and 
\[
\kappa_n = \sum_{\substack{\ell=0 \\ \ell \neq n}}^3    \chi_{\ell,0}(\sigma_\ell , x_\ell)
.\]

{Therefore,  by adding and  subtracting the same terms and ordering them appropriately we   have the following:}
\begin{align*}
{}^{\psi}\mathcal D^{\sigma,\beta}_{\nu} [f] (x)
=&   \sum_{n=0 }^3 \psi_n     \frac{\partial }{\partial x_n} \left[ \varsigma(x)   \left( \sum_{k=0 }^3   \lambda_{k}(f ) (x)  \right) \right]    -
 \sum_{\substack{n, k = 0 \\ n \neq k}}^3 \psi_n     \frac{\partial }{\partial x_n} \left[ \varsigma(x)       \lambda_{k}(f ) (x) \right]\\ 
 &-    \sum_{n=0 }^3 \psi_n     \frac{\partial }{\partial x_n} \left[ \sum_{k=0}^3 \kappa_k     \lambda_{k}(f ) (x) \right] 
  + 
   \sum_{\substack{n, k = 0 \\ n\neq k }}^3     \psi_n   \frac{\partial }{\partial x_n} \left[  \kappa_k      \lambda_{k}(f ) (x) \right]\\ 
  &- \sum_{{ n, m = 0} }^3 \psi_n\psi_m  L_{n,m}[f](x) \lambda_{n}(f_m) (x) , 
	\end{align*}
for all ${x\in\Omega}$. As a consequence we have that 
 \begin{align*}
{}^{\psi}\mathcal D^{\sigma,\beta}_{\nu} [f] (x)
 = { }&{ }  \sum_{n=0 }^3 \psi_n     \frac{\partial }{\partial x_n} \left[ \varsigma(x)   \left( \sum_{k=0 }^3   \lambda_{k}(f ) (x)  \right) - \sum_{k=0}^3 \kappa_k      \lambda_{k}(f ) (x) \right]  \\
 { }&{ } + \sum_{\substack{ n, k = 0 \\  n \neq k}}^3 \psi_n     \frac{\partial }{\partial x_n} \left[ \kappa_k      \lambda_{k}(f ) (x)-  \varsigma(x)       \lambda_{k}(f ) (x) \right] \\
  { }&{ } - \sum_{n, m = 0}^3 \psi_n\psi_m  L_{n,m}[f](x) \lambda_{n}(f_m) (x),
\end{align*} 
i.e.,
 \begin{align*}
{}^{\psi}\mathcal D^{\sigma,\beta}_{\nu} [f] (x)
 = { }&{ } \sum_{n=0 }^3 \psi_n     \frac{\partial }{\partial x_n} \left[\sum_{k=0 }^3 (\chi_{k,0}(\sigma_k, x_k )) \lambda_{k}(f) (x) \right]
\\
{ }&{ }- \sum_{\substack{n, k = 0 \\ n\neq k}}^3 \psi_n \frac{\partial }{\partial x_n} \left[ ( \chi_{k,0}(\sigma_k, x_k )   )     \lambda_{k}(f ) (x)  \right]\\
{ }&{ }-\sum_{n, m = 0}^3 \psi_n\psi_m  L_{n,m}[f](x) \lambda_{n}(f_m)(x), 
\end{align*}
or equivalently 
 \begin{align*}
{}^{\psi}\mathcal D^{\sigma,\beta}_{\nu} [f] (x)  = { }&{ }  {}^{\psi} \mathcal D \circ \mathfrak L ^{\sigma,\beta}_{\nu}[f](x) -
\sum_{\substack{n, k = 0 \\ n\neq k}}^3 \psi_n \frac{\partial }{\partial x_n} \left[ (\chi_{k,0}(\sigma_k, x_k ))  \lambda_{k}(f)(x) \right] 
 \\
{ }&{ }-\sum_{n, m = 0}^3 \psi_n\psi_m  L_{n,m}[f](x) \lambda_{n}(f_m) (x),
\end{align*}
for all ${x\in \Omega}$.
\end{proof}

Under the hypothesis and notation of the previous proposition, let us present some consequences of quaternionic Borel--Pompeiu and Stokes formulas
given in \eqref{BorelHyp} and \eqref{StokesHyp}. 
\begin{proposition}\label{propBPS}
Let ${\Omega\subset\mathbb H}$ be a domain such that ${\partial \Omega}$ is a 3-dimensional compact smooth surface. If 

\[
\mathfrak L ^{\sigma,\beta}_{\nu}  [f] ,  \mathfrak L ^{\rho,\delta}_{\mu}[g]  \in C^1(\Omega,\mathbb H)
,\]
then 
 \begin{align}\label{BP-New}   \int_{\partial \Omega}&(K_{\psi}(\tau-x)\sigma_{\tau}^{\psi}  \mathfrak L ^{\sigma,\beta}_{\nu}  [f] (\tau)  +   \mathfrak L ^{\rho,\delta}_{\mu}[g](\tau)   \sigma_{\tau}^{\psi} K_{\psi}(\tau-x) ) \nonumber  \\ 
&  - 
{
\int_{\Omega} \left(K_{\psi} (y-x)    
 {}^{\psi}\mathcal D^{\sigma,\beta}_{\nu} [f] (y)   + 
  {}^{\psi}\mathcal D^{\rho,\delta}_{r,\mu} [g] (y)    K_{\psi} (y-x)
    \right)  dy   \nonumber   }\\
    &  -   
\int_{\Omega}  K_{\psi} (y-x)    \left(
  \mathcal E^{\sigma,\beta}_{\nu}  [f](y)    +
 \sum_{{ n, m = 0} }^3 \psi_n\psi_m  L_{n,m}[f](y) \lambda^{\beta_n}_{\nu_n}(f_m) (y) \right) dy   \nonumber   \\
   &   - 
\int_{\Omega} \left( 
\mathcal E^{\rho,\delta}_{r,\mu}  [g](y) +
 \sum_{{ n, m = 0} }^3 \psi_m\psi_n  T_{n,m}[g](y) \lambda^{\delta_n}_{\mu_n}(g_m) (y) \right)   K_{\psi} (y-x)
      dy   \nonumber \\
		=  &  \left\{ \begin{array}{ll}   \mathfrak L ^{\sigma,\beta}_{\nu}  [f] (x) +  \mathfrak L ^{\rho,\delta}_{\mu}[g] (x) , &  x\in \Omega,  \\ 0 , &  x\in \mathbb H\setminus\overline{\Omega},                     
\end{array} \right. 
\end{align} 
and 
\begin{align}\label{New-Stokes} 
\int_{\partial \Omega}  \mathfrak L ^{\rho,\delta}_{\mu}[g]  \sigma^\psi_x & \mathfrak L ^{\sigma,\beta}_{\nu}  [f] \nonumber \\
= { }&{ } \int_{\Omega } \left(  \mathfrak L ^{\rho,\delta}_{\mu}[g]   
		 		  {}^{\psi}\mathcal D^{\sigma,\beta}_{\nu} [f]    +  {}^{\psi}\mathcal D^{\rho,\delta}_{r,\mu} [g]  
		 		    		  \mathfrak L ^{\sigma,\beta}_{\nu}  [f]   \right)dx +
  \nonumber   \\
  &{ } +  \int_{\Omega }    \mathfrak L ^{\rho,\delta}_{\mu}[g]  \left(\mathcal E^{\sigma,\beta}_{\nu}  [f]  
  +  \sum_{{ n, m = 0} }^3 \psi_n\psi_m  L_{n,m}[f]  \lambda^{\beta_n}_{\nu_n}(f_m)  
			  \right)dx   \nonumber   \\ 
	&{ } + \int_{\Omega } \left(  \mathcal E^{\rho,\delta}_{r,\mu}  [g]  +
 \sum_{{ n, m = 0} }^3 \psi_m\psi_n  T_{n,m}[g]  \lambda^{\delta_n}_{\mu_n}(g_m) 
		 \right)		  \mathfrak L ^{\sigma,\beta}_{\nu}  [f]  dx.
\end{align}
\end{proposition}
\begin{proof}
The formulas follow by applying the quaternionic Borel--Pompeiu and Stokes formulas to the functions 
${\mathfrak L ^{\sigma,\beta}_{\nu} [f]}$,  ${\mathfrak L ^{\rho,\delta}_{\mu}[g]}$ and using the identities in  Proposition \ref{Fueter-Hdefpropbfractal}.  
\end{proof}

\begin{remark}\label{REMCAU}
If  ${\mathfrak L ^{\sigma,\beta}_{\nu}}$ and   ${\mathfrak L ^{\rho,\delta}_{\mu}}$  are invertible operators we can use formula \eqref{BP-New} to obtain the quaternionic values of ${f}$ and  ${g}$. In addition, if 
\[
f\in   {}^{\psi}\mathcal  M^{\sigma,\beta}_{\nu} (\Omega,\mathbb H),
\]
 and 
\[
g\in   {}^{\psi}\mathcal M^{\rho,\delta}_{r,\mu} (\Omega,\mathbb H)
,\]
 then 
\begin{align*}  &  \int_{\partial \Omega}(K_{\psi}(\tau-x)\sigma_{\tau}^{\psi}  \mathfrak L ^{\sigma,\beta}_{\nu}  [f] (\tau)  +   \mathfrak L ^{\rho,\delta}_{\mu}[g](\tau)   \sigma_{\tau}^{\psi} K_{\psi}(\tau-x) ) \nonumber  \\ 
    &  -   
\int_{\Omega}  K_{\psi} (y-x)    \left(
  \mathcal E^{\sigma,\beta}_{\nu}  [f](y)    +
 \sum_{{ n, m = 0} }^3 \psi_n\psi_m  L_{n,m}[f](y) \lambda^{\beta_n}_{\nu_n}(f_m) (y) \right) dy   \nonumber   \\
   &   - 
\int_{\Omega} \left( 
\mathcal E^{\rho,\delta}_{r,\mu}  [g](y) +
 \sum_{{ n, m = 0} }^3 \psi_m\psi_n  T_{n,m}[g](y) \lambda^{\delta_n}_{\mu_n}(g_m) (y) \right)   K_{\psi} (y-x)
      dy   \nonumber \\
		=  &  \left\{ \begin{array}{ll}   \mathfrak L ^{\sigma,\beta}_{\nu}  [f] (x) +  \mathfrak L ^{\rho,\delta}_{\mu}[g] (x) , &  x\in \Omega,  \\ 0 , &  x\in \mathbb H\setminus\overline{\Omega},                     
\end{array} \right. 
\end{align*} 
and 
\begin{align*} 
&		\int_{\partial \Omega}  \mathfrak L ^{\rho,\delta}_{\mu}[g]  \sigma^\psi_x  \mathfrak L ^{\sigma,\beta}_{\nu}  [f]   =    \int_{\Omega }    \mathfrak L ^{\rho,\delta}_{\mu}[g]  \left(\mathcal E^{\sigma,\beta}_{\nu}  [f]  
  +  \sum_{{ n, m = 0} }^3 \psi_n\psi_m  L_{n,m}[f]  \lambda^{\beta_n}_{\nu_n}(f_m)  
			  \right)dx   \nonumber   \\ 
	& + \int_{\Omega } \left(  \mathcal E^{\rho,\delta}_{r,\mu}  [g]  +
 \sum_{{ n, m = 0} }^3 \psi_m\psi_n  T_{n,m}[g]  \lambda^{\delta_n}_{\mu_n}(g_m) 
		 \right)		  \mathfrak L ^{\sigma,\beta}_{\nu}  [f]  dx.
\end{align*}
\end{remark}
\begin{remark}
For some specific functions ${\sigma}$, ${\beta}$ and ${\nu}$ we can see that ${{}^{\psi}\mathcal D^{\sigma,\beta}_{\nu}}$  coincides with ${{}^{\psi}\mathcal D}$. Therefore,     Proposition \ref{propBPS} and Remark \ref{REMCAU} can give us  the  usual Borel--Pompeiu, Stokes and Cauchy formulas presented   in \cite{SV1,SV2}.
\end{remark}

\section{Quaternionic  proportional fractal truncated Fueter operator}\label{sec:truncated}
We now consider the partial differential operators given by Remarks \ref{rem1} and \ref{rem2}. Let ${{k}:=(k_0,k_1,k_2,k_3)\in \mathbb N^4}$ (or ${k=(\infty, \infty, \infty, \infty)}$), 
\[
\sigma = (\sigma_0,\sigma_1,\sigma_2,\sigma_3)\in (0,1]^4
,\]
 {
\[
\beta =(\beta_0,\beta_1, \beta_2,\beta_3)\in (0,1]^4
,\]
}  
\[
\alpha =(\alpha_0,\alpha_1, \alpha_2,\alpha_3)\in (0,1]^{4},
\]
  with coordinate index ${n=0,1,2,3}$.  

Recall that 
\[
\mathbb H^+:=\{x= x_0+ x_1e_1+ x_2e_2 + x_3e_3 \in \mathbb H \ \mid \  x_0,x_1,x_2,x_3 >0\}
.\]
 We consider domains ${\Omega\subset \mathbb H}$ and functions ${f,g :\Omega\to \mathbb H}$ satisfying ${f(\Omega),g(\Omega) \subset \mathbb H^{+}}$. The proportional ${\beta_n}$-fractal partial derivative of ${f\in C^{1}(\Omega,\mathbb{R})}$ is given by
\[\dfrac{ {\partial }^{\sigma_n,\beta_n}f}{ \partial x _{\alpha_n, k_n} } (x) : = (1-\sigma_n)f (x) +   \sigma_n  \frac{ \dfrac{\partial  f ^{\beta_n} }{\partial x_n} (x) }{ \dfrac{ \partial  e ( x_n^{\alpha_n})_{k_n } }{\partial x_n  } },\]
for all ${x=\sum_{n=0}^3 \psi_{n} x_n \in \Omega}$. 

In addition, let ${{m}:=(m_0,m_1,m_2,m_3)\in \mathbb N^4}$ (or ${m=(\infty, \infty, \infty, \infty)}$), 
\[
\rho = (\rho _0,\rho _1,\rho _2,\rho _3) \in (0,1]^{4}
,\]
 
\[
\delta =(\delta_0,\delta_1, \delta_2,\delta_3)\in (0,1]^4
,\]
 
\[
\gamma =(\gamma_0,\gamma_1, \gamma_2,\gamma_3) \in (0,1]^{4}
.\]
 For each coordinate index ${n=0,1,2,3}$, the proportional ${\delta_n}$-fractal partial derivative of ${g\in C^1(\Omega,\mathbb R)}$ is given by 
\[\dfrac{ {\partial }^{\rho_n,\delta_n}g}{ \partial x _{\gamma_n, m_n} } (x) : = (1-\rho_n) g(x) +   \rho_n \frac{ \dfrac{\partial  g ^{\delta_n} }{\partial x_n} (x) }{ \dfrac{ \partial  e ( x_n^{\gamma_n})_{m_n } }{\partial x_n  } },\]
for all ${x=\sum_{n=0}^3 \psi_{n} x_n \in \Omega}$. 
We will define the  ${\psi}$-proportional-fractal Fueter-type operators associated with   truncated exponential functions. 
\begin{definition}\label{def100}  
Let 
\[
f= \sum_{\ell=0}^3 \psi_\ell  f_\ell, \ g= \sum_{\ell=0}^3 \psi_\ell  g_\ell \in C^1(\Omega,\mathbb H)
,\]
 where ${f_0,f_1,f_2,f_3,g_0,g_1,g_2,g_3}$ are ${\mathbb{R}}$-valued functions. Define the ${\psi}$-proportional-fractal Fueter-type operators associated with truncated exponential functions as follows: 
\begin{align}\label{PFFl} 
{}^{\psi}{\mathcal D}^{\sigma,\beta}_{\alpha, k}[f](x):= & \sum_{n, \ell = 0 }^3 \psi_n \psi_\ell \frac{{\partial }^{\sigma_n,\beta_n} f_\ell}{\partial x_{\alpha_n, k_n}}(x),
\end{align}
and 
\begin{align*} 
{}^{\psi}{\mathcal D}^{\rho,\delta}_{r, \gamma, m}[g](x) := & \sum_{n, \ell = 0 }^3 \psi_\ell \frac{{\partial }^{\rho_n,\delta_n} g_\ell}{\partial x_{\gamma_n, m_n}}(x) \psi_n.
\end{align*} 

To simplify notation, we introduce the following auxiliary operators:   
\begin{align*}
H^{\sigma_n, \beta_n }_{\alpha_n, k_n}[f_\ell](x )  :=  & \int_{0}^{x_n} \frac{1- \sigma_{n} }{\sigma_n} \left( \frac{d }{dx_n} e ( x_n^{\alpha_n})_{k _n} \right) f_\ell (x)^{ 1-\beta_n } dx_n,  \\
H^{\rho_n, \delta_n }_{\gamma_n, m_n}[g_{\ell}](x )  :=  & \int_{0}^{x_n} \frac{1-\rho_{n} }{\rho_n} \left( \frac{d }{dx_n} e ( x_n^{\gamma_n})_{m _n} \right) g_\ell (x)^{ 1-\delta_n } dx_n . 
\end{align*}
Use 
\[
h_{n,\ell}(x)= H^{\sigma_n, \beta_n }_{\alpha_n, k_n}(f_\ell)(x),
\]
 and 
\[
j_{n,\ell}(x) = H^{\rho_n, \delta_n}_{\gamma_n, m_n}(g_\ell)(x),
\]
 for all ${n,\ell =0,1,2,3}$. Denote 
\begin{align*}
T^{\sigma_n, \beta_n }_{\alpha_n, k_n} [f_\ell] (x):= & \frac{\partial }{\partial x_n} \left(  \frac{\sigma_n}{e^{h_{n ,\ell}(x)}    
\dfrac{d }{dx_n} e ( x_n^{\alpha_n})_{k _n}} \right) e^{h_{n,\ell}(x)}f_\ell (x)^{\beta_n},\\
{}^{\psi}W_{\alpha, k}^{\sigma, \beta} [ f ] (x)   
 :=  &
\sum_{{
 \substack{\ell, n = 0 \\ \ell \ne n}} }^3 \frac{\sigma_\ell}{\dfrac{d }{dx_\ell} e(x_\ell^{\alpha_\ell})_{k _\ell}} \psi_n \frac{\partial}{\partial x_n}I^{\beta_\ell} [f](x), \\
S^{\rho_n, \delta_n }_{\gamma_n, m_n} [g_\ell] (x):= & \frac{\partial }{\partial x_n} \left(  \frac{\rho_n}{e^{j_{n ,\ell}(x)   }    
   \dfrac{d }{dx_n} e ( x_n^{\gamma_n})_{m _n}}\right) e^{j_{n,\ell}(x)} g_\ell (x)^{\delta_n},\\
 {}^{\psi}V_{\gamma, m}^{\rho, \delta} [g] (x):=  & \sum_{
 {
 \substack{\ell, n = 0 \\ \ell \ne n}}}^3 \frac{\rho_\ell}{     
\dfrac{d }{dx_\ell} e (x_\ell^{\gamma_\ell})_{m _\ell}}\frac{\partial }{\partial x_n} I^{\delta_\ell} [g](x) \psi_n,
\end{align*}
where 
\[
\displaystyle I^{\beta_n}[f] (x) = \sum_{ \ell=0 }^3 \psi_\ell f_\ell  (x)^{\beta_n},
\]
 for all ${x \in \Omega}$ and  ${n=0,1,2,3}$. The operator ${I^{\delta_n}[g](x)}$ is defined analogously. 
\end{definition}

{It should be noted that  our scalars are quaternions and these form a skew field. Thus, the function sets under consideration are quaternionic modules, with a distinction between right and left modules because quaternionic multiplication is noncommutative.}

\begin{definition}
For ${\beta=(1,1,1,1)}$, we see that  the quaternionic-right module, or ${\mathbb H}$-right module, of ${\psi}$-proportional-fractal hyperholomorphic functions associated with truncated exponential functions  is given by 
 
\[
{}^{\psi}{\mathcal M}^{\sigma,\beta}_{\alpha, k}(\Omega,\mathbb H)  
  = C^1(\Omega, \mathbb H)  \cap \textrm{Ker} {}^{\psi}{\mathcal D}^{\sigma,\beta}_{\alpha, k}
.\]
 Meanwhile, the ${\mathbb H}$-left module of  right   ${\psi}$-proportional-fractal hyperholomorphic functions   associated with   truncated exponential functions  is given by 

\[
{}^{\psi}{\mathcal M}^{\rho,\delta}_{r, \gamma, m} (\Omega, \mathbb H) =
C^1  (\Omega, \mathbb H) \cap \textrm{Ker} {}^{\psi}{\mathcal D}^{\rho,\delta}_{r, \gamma, m}
.\]

\end{definition}
If ${\beta=(1,1,1,1)}$, then the condition ${f(\Omega) \subset \mathbb H^{+}}$
is not necessary.

\begin{remark}
{At a very superficial level, we can observe that for certain expressions  ${\alpha}$, ${\beta}$, ${\sigma}$ and ${k}$, we obtain ${{}^{\psi}{\mathcal D}}$ from 
\[
{}^{\psi}{\mathcal D}^{\sigma,\beta}_{\alpha, k},
\]
 as a special case. Then  we obtain the well-known formulas 
   \[\overline{{}^{\psi}{\mathcal D}} \circ {}^{\psi}{\mathcal D}=  
   {}^{\psi}{\mathcal D}\circ \overline{{}^{\psi}{\mathcal D}} = \Delta_{\mathbb R^4},\] 
  where  ${\Delta_{\mathbb R^4}}$ is the four-dimensional Laplace operator; see \cite{SV1}. These formulas are obtained as special cases of the compositions:
    \[ \overline{{}^{\psi}{\mathcal D}^{\sigma,\beta}_{\alpha, k} }\circ {}^{\psi}{\mathcal D}^{\sigma,\beta}_{\alpha, k}, \quad   {}^{\psi}{\mathcal D}^{\sigma,\beta}_{\alpha, k} \circ \overline{{}^{\psi}{\mathcal D}^{\sigma,\beta}_{\alpha, k} }.\]}

On the other hand, note that   
\[
f= \sum_{\ell=0}^3 \psi_\ell  f_\ell \in   {}^{\psi}{\mathcal M}^{\sigma,\beta}_{\alpha, k}(\Omega,\mathbb H),
\]
 if and only if the real components of ${f}$ satisfy the Cauchy--Riemann type equation extending the  Cauchy--Riemann type  equations given in \cite{SV1}.
For example, if ${\beta =(1,1,1,1)}$ and ${k=(1,1,1,1)}$, then  Cauchy--Riemann type  equations  are the following:  
     \begin{align*} 
  0=& (1-\sigma_0)f_{0} (x) +   \sigma_0 \frac{ \dfrac{\partial   f_0     }{\partial x_0}  (x)  }{  {\alpha_0}  x_0^{\alpha_0-1}    }
   -       (1-\sigma_1)f_{1} (x) -   \sigma_1 \frac{ \dfrac{\partial   f_1     }{\partial x_1}  (x)  }{  {\alpha_1}  x_1^{\alpha_1-1}    }   \\
& -        (1-\sigma_2)f_{2} (x) -  \sigma_2 \frac{ \dfrac{\partial   f_2     }{\partial x_2}  (x)  }{  {\alpha_2}  x_2^{\alpha_2-1}    }
          -            (1-\sigma_3)f_{3} (x) -   \sigma_3 \frac{ \dfrac{\partial   f_3     }{\partial x_3}  (x)  }{  {\alpha_3}  x_3^{\alpha_3-1}},  
			\\
         0=& (1-\sigma_0)f_{1} (x) +   \sigma_0 \frac{ \dfrac{\partial   f_1     }{\partial x_0}  (x)  }{  {\alpha_0}  x_0^{\alpha_0-1}    }
   +     (1-\sigma_1)f_{0} (x) +  \sigma_1 \frac{ \dfrac{\partial   f_0     }{\partial x_1}  (x)  }{  {\alpha_1}  x_1^{\alpha_1-1}    }   \\
& +    (1-\sigma_2)f_{3} (x) +  \sigma_2 \frac{ \dfrac{\partial   f_3     }{\partial x_2}  (x)  }{  {\alpha_2}  x_2^{\alpha_2-1}    }
         -            (1-\sigma_3)f_{2} (x) -   \sigma_3 \frac{ \dfrac{\partial   f_2     }{\partial x_3}  (x)  }{  {\alpha_3}  x_3^{\alpha_3-1}    }
       ,  
			\\
                0=& (1-\sigma_0)f_{2} (x) +   \sigma_0 \frac{ \dfrac{\partial   f_2     }{\partial x_0}  (x)  }{  {\alpha_0}  x_0^{\alpha_0-1}    }
   +     (1-\sigma_2)f_{0} (x) +  \sigma_2 \frac{ \dfrac{\partial   f_0     }{\partial x_2}  (x)  }{  {\alpha_2}  x_2^{\alpha_2-1}    }   \\
&+    (1-\sigma_3)f_{1} (x) +  \sigma_3 \frac{ \dfrac{\partial   f_1     }{\partial x_3}  (x)  }{  {\alpha_3}  x_3^{\alpha_3-1}    } 
         -            (1-\sigma_1)f_{3} (x) -   \sigma_1 \frac{ \dfrac{\partial   f_3     }{\partial x_1}  (x)  }{  {\alpha_1}  x_1^{\alpha_1-1}    }
       ,        \\
       0=& (1-\sigma_0)f_{3} (x) +   \sigma_0 \frac{ \dfrac{\partial   f_3     }{\partial x_0}  (x)  }{  {\alpha_0}  x_0^{\alpha_0-1}    }
   +     (1-\sigma_3)f_{0} (x) +  \sigma_3 \frac{ \dfrac{\partial   f_0     }{\partial x_3}  (x)  }{  {\alpha_3}  x_3^{\alpha_3-1}    }   \\
& 
+ (1-\sigma_1)f_{2} (x) + \sigma_1 \frac{ \dfrac{\partial f_2}{\partial x_1} (x)}{{\alpha_1} x_1^{\alpha_1-1}}
 - (1-\sigma_2)f_{1} (x) - \sigma_2 \frac{ \dfrac{\partial f_1}{\partial x_2} (x)}{{\alpha_2} x_2^{\alpha_2-1}},
\end{align*}
for all ${x \in \Omega}$.
\end{remark}

\begin{example}\label{example1} 
If ${\beta=(1,1,1,1)}$ and ${k=(1,1,1,1)}$, then \eqref{PFFl} is
\begin{align*} 
{}^{\psi}{\mathcal D}^{\sigma,1}_{\alpha, 1}[f](x)= & \sum_{n, \ell = 0 }^3 \psi_n \psi_\ell \frac{{\partial }^{\sigma_n,1} f_\ell}{\partial x_{\alpha_n, 1}}(x)=\sum_{n, \ell = 0 }^3 \psi_n \psi_\ell\left( (1-\sigma_n)f_{\ell} (x) +   \sigma_n  \frac{ \dfrac{\partial  f_{\ell}   }{\partial x_n} (x) }{ \dfrac{ \partial  e ( x_n^{\alpha_n})_{1} }{\partial x_n  } }\right).
\end{align*}
Define the function ${f(x)= \frac{5}{2} x + \frac{3}{2}\bar x}$, for all 
\[
x=\sum_{\ell=0}^3\psi_\ell x_\ell \in \mathbb H,
\]
 such that {${x_0,x_1,x_2,x_3> 0}$}.  Then we see that {${f=\sum_{\ell =0}^3 \psi_\ell \zeta_\ell}$, where
\begin{equation*}
    \zeta_{0}(x)=x_{0},\quad \zeta_{1}(x)=x_{1}-\psi_{1}x_{0},\quad \zeta_{2}(x)=x_{2}-\psi_{2}x_{0},\quad \zeta_{3}(x)=x_{3}-\psi_{3}x_{0}.  
\end{equation*}
The functions  ${\zeta_1, \ \zeta_2, \ \zeta_3}$} are known as Fueter variables and are used to decompose the hyperholomorphic functions in terms of series of symmetric polynomials. 
 Note that {
\[
f(x)=\frac{5}{2} x + \frac{3}{2}\bar x =  4x_0 +\psi_1x_1 + \psi_2x_2+ \psi_3x_3,
\]
  for all ${x}$; i.e., ${f_0(x)=4x_0}$ and ${f_\ell (x) = x_\ell}$ for ${\ell=1,2,3}$} and
\begin{align*}
    {}^{\psi}{\mathcal D}^{\sigma,1}_{\alpha, 1}[f](x)=&(1-\sigma_0)f   (x) 
     +   \psi_1  (1-\sigma_1)f (x) + \psi_2(1-\sigma_2)f (x)  + \psi_3    (1-\sigma_3)f (x) \\
    &     +   \sigma_0 \frac{ \dfrac{\partial   f     }{\partial x_0}  (x)  }{  {\alpha_0}  x_0^{\alpha_0-1}    }
  + \psi_1 \sigma_1 \frac{ \dfrac{\partial   f     }{\partial x_1}  (x)  }{  {\alpha_1}  x_1^{\alpha_1-1}    } 
   +   \psi_2\sigma_2 \frac{ \dfrac{\partial   f     }{\partial x_2}  (x)  }{  {\alpha_2}  x_2^{\alpha_2-1}    }
  +  \psi_3\sigma_3 \frac{ \dfrac{\partial   f     }{\partial x_3}  (x)  }{  {\alpha_3}  x_3^{\alpha_3-1}    }\\
  =&\bigg( (1-\sigma_0)      +   \psi_1  (1-\sigma_1)  +\psi_2 (1-\sigma_2)   + \psi_3    (1-\sigma_3)\bigg) \bigg( \frac{5}{2} x + \frac{3}{2}\bar x  \bigg)  \\
    &     +   \sigma_0 \frac{4}{  {\alpha_0}  x_0^{\alpha_0-1}    }
  + \psi_1 \sigma_1 \frac{ \psi_1 }{  {\alpha_1}  x_1^{\alpha_1-1}    } 
   +   \psi_2\sigma_2 \frac{ \psi_2 }{  {\alpha_2}  x_2^{\alpha_2-1}    }
  +  \psi_3\sigma_3 \frac{ \psi_3 }{  {\alpha_3}  x_3^{\alpha_3-1}    }.\end{align*}
Therefore, 
\begin{align*}
    {}^{\psi}{\mathcal D}^{\sigma,1}_{\alpha, 1}[f](x)=&   \bigg( (1-\sigma_0)      +   \psi_1  (1-\sigma_1)  +\psi_2 (1-\sigma_2)   + \psi_3    (1-\sigma_3)\bigg) \bigg( \frac{5}{2} x + \frac{3}{2}\bar x  \bigg)  \\
    &     +   \sigma_0 \frac{4}{  {\alpha_0}  x_0^{\alpha_0-1}    }
  - \sigma_1 \frac{ 1}{  {\alpha_1}  x_1^{\alpha_1-1}    } 
   -   \sigma_2 \frac{ 1}{  {\alpha_2}  x_2^{\alpha_2-1}    }
  -\sigma_3 \frac{ 1 }{  {\alpha_3}  x_3^{\alpha_3-1}    },
\end{align*}
 for all 
\[
x=\sum_{\ell=0}^3\psi_\ell x_\ell \in \mathbb H,
\]
 such that { ${x_0,x_1,x_2,x_3> 0}$}.
 \end{example}

\begin{proposition}\label{prop1234} 
Suppose that  
\[
\sigma = (\sigma_0,\sigma_1,\sigma_2,\sigma_3), \ \rho = (\rho_0,\rho_1,\rho_2,\rho_3)\in (0,1]^4
.\]
Let 
\[
f= \sum_{\ell=0}^3 \psi_\ell  f_\ell, \ g= \sum_{\ell=0}^3 \psi_\ell  g_\ell  \in C^1(\Omega,\mathbb H)
.\]
 Then  
\begin{align}\label{equa10} 
& {}^{\psi}\mathcal D \left[ \sum_{\ell=0}^3  \frac{\sigma_\ell}{     
   \dfrac{d }{dx_\ell} e ( x_\ell^{\alpha_\ell})_{k _\ell}    }    I^{\beta_\ell} [f]     \right] (x) \nonumber  \\ 
   =  & {}^{\psi}{\mathcal D}^{\sigma,\beta}_{\alpha, k}[f](x)  +   \sum_{n, \ell = 0 }^3  \psi_n \psi_\ell 
    T^{\sigma_n, \beta_n }_{\alpha_n, k_n} [f_\ell] (x)  + 
{}^{\psi}W_{\alpha, k}^{\sigma, \beta} [ f ] (x),   
\end{align}
and 
\begin{align}\label{equa11} 
& {}^{\psi}\mathcal D_r \left[ \sum_{\ell=0}^3  \frac{\rho_\ell}{     
   \dfrac{d }{dx_\ell} e ( x_\ell^{\gamma_\ell})_{m _\ell}    }    I^{\delta_\ell} [g]     \right] (x) \nonumber  \\ 
   =  & {}^{\psi}{\mathcal D}^{\rho,\delta}_{r, \gamma, m}[g](x)  +   \sum_{n, \ell = 0 }^3   \psi_\ell 
    S^{\rho_n, \delta_n }_{\gamma_n, m_n} [g_\ell] (x)  \psi_n + 
{}^{\psi}V_{\gamma, m}^{\rho, \delta} [g ] (x),   \end{align}
for all ${x \in \Omega}$.
\end{proposition}

\begin{proof} 
The proofs of the two identities are similar, so we prove only the first. We have
\begin{align*} \frac{\partial }{\partial x_n} \left(   e^{h_{n,\ell}(x_n)}   f_\ell  (x)^{\beta_n} \right)   = &    \left[\  \frac{1-\sigma_n}{\sigma_n}  \left( \frac{d }{dx_n} e ( x_n^{\alpha_n})_{k _n} \right)  (f_\ell (x) )^{1- \beta_n}  f_{\ell} ^{\beta_n}(x)  
 + \frac{\partial  f_\ell ^{\beta_n} }{\partial x_n}  (x)    \right] e^{h_{n,\ell }(x )} \\
 =  &   \left[ (1-\sigma_n)f_{\ell} (x) +   \sigma_n \frac{ \dfrac{\partial f_\ell^{\beta_n}}{\partial x_n}(x)}{\dfrac{d }{dx_n} e(x_n^{\alpha_n})_{k _n}} \right] \frac{1}{\sigma_n} e^{h_{n,\ell}(x)}  
\frac{d }{dx_n} e ( x_n^{\alpha_n})_{k _n}, \\
\frac{ {\partial }^{\sigma_n,\beta_n}}{ \partial x_{\alpha_n, k_n} }f_\ell (x) = &
\frac{\sigma_n}{e^{h_{n,\ell}(x)} \dfrac{d}{dx_n} e( x_n^{\alpha_n})_{k _n}} \frac{\partial}{\partial x_n} \left(e^{h_{n,\ell}(x)} 
f_\ell (x)^{\beta_n} \right),   
\end{align*} 
and
\begin{align*} 
{}^{\psi}{\mathcal D}^{\sigma,\beta}_{\alpha, k} [f](x) = & \sum_{n=0}^3 \psi_n\frac{{\partial }^{\sigma_n,\beta_n}}{\partial x _{\alpha_n, k_n} }f (x) \\
= & \sum_{n, \ell = 0 }^3  \psi_n \psi_\ell \frac{\sigma_n}{e^{h_{n,\ell}(x)} \dfrac{d}{dx_n} e( x_n^{\alpha_n})_{k _n}} \frac{\partial}{\partial x_n} \left(e^{h_{n,\ell}(x )} f_\ell (x)^{\beta_n} \right).
\end{align*}
	
The identities  
\begin{align*}
 &  \frac{\partial }{\partial x_n} \left(  \frac{\sigma_n}{     
   \dfrac{d }{dx_n} e ( x_n^{\alpha_n})_{k _n}    }        f_\ell  (x)^{\beta_n} \right) =  
   \frac{\partial }{\partial x_n} \left(  \frac{\sigma_n}{e^{h_{n,\ell }(x)   }    
   \dfrac{d }{dx_n} e ( x_n^{\alpha_n})_{k _n}    }       e^{h_{n,\ell }(x )}   f_\ell  (x)^{\beta_n} \right)  \\
 = &  \frac{\partial }{\partial x_n} \left(  \frac{\sigma_n}{e^{h_{n,\ell }(x)   }    
   \dfrac{d }{dx_n} e ( x_n^{\alpha_n})_{k _n}    }  \right)     e^{h_{n,\ell}(x )}   f_\ell  (x)^{\beta_n}    +
  \frac{\sigma_n}{e^{h_{n ,\ell}(x)   }    
   \dfrac{d }{dx_n} e ( x_n^{\alpha_n})_{k _n}    }   \frac{\partial }{\partial x_n} \left(     e^{h_{n,\ell}(x )}   f_\ell  (x)^{\beta_n} \right)  ,
 \end{align*}
 and 
\begin{align*}
& \frac{\sigma_n}{e^{h_{n,\ell }(x)   }    
   \dfrac{d }{dx_n} e ( x_n^{\alpha_n})_{k _n}    }   \frac{\partial }{\partial x_n} \left(     e^{h_{n,\ell}(x )}   f_\ell  (x)^{\beta_n} \right)  \\
     =&  \frac{\partial }{\partial x_n} \left(  \frac{\sigma_n}{     
   \dfrac{d }{dx_n} e ( x_n^{\alpha_n})_{k _n}    }        f_\ell  (x)^{\beta_n} \right) -  
     \frac{\partial }{\partial x_n} \left(  \frac{\sigma_n}{e^{h_{n ,\ell}(x)   }    
   \dfrac{d }{dx_n} e ( x_n^{\alpha_n})_{k _n}    }  \right)     e^{h_{n,\ell}(x)}   f_\ell  (x)^{\beta_n}   \\
    =&  \frac{\partial }{\partial x_n} \left(  \frac{\sigma_n}{     
   \dfrac{d }{dx_n} e ( x_n^{\alpha_n})_{k _n}    }        f_\ell  (x) ^{\beta_n} \right) -  
   T^{\sigma_n, \beta_n }_{\alpha_n, k_n} [f_\ell] (x)   (x)  ,   
 \end{align*}
imply that 
\begin{align*} 
{}^{\psi}{\mathcal D}^{\sigma,\beta}_{\alpha, k} [f](x)  = &  \sum_{n, \ell = 0 }^3  \psi_n \psi_\ell 
\left[  \frac{\partial }{\partial x_n} \left(  \frac{\sigma_n}{     
   \dfrac{d }{dx_n} e ( x_n^{\alpha_n})_{k _n}    }        f_\ell  (x)^{\beta_n} \right) -  
    T^{\sigma_n, \beta_n }_{\alpha_n, k_n} [f_\ell]    (x)    \right]
    \\
     =&  \sum_{n, \ell = 0 }^3  \psi_n \psi_\ell    \frac{\partial }{\partial x_n} \left[  \frac{\sigma_n}{     
   \dfrac{d }{dx_n} e ( x_n^{\alpha_n})_{k _n}    }        f_\ell  (x)^{\beta_n} \right]  -   \sum_{n, \ell = 0 }^3  \psi_n \psi_\ell 
    T^{\sigma_n, \beta_n }_{\alpha_n, k_n} [f_\ell]   (x) \\
         =&  \sum_{n=0 }^3  \psi_n     \frac{\partial }{\partial x_n} \left[  \frac{\sigma_n}{     
   \dfrac{d }{dx_n} e ( x_n^{\alpha_n})_{k _n}    }   \sum_{ \ell=0 }^3 \psi_\ell       f_\ell  (x)^{\beta_n} \right] -   \sum_{n, \ell = 0 }^3  \psi_n \psi_\ell 
    T^{\sigma_n, \beta_n }_{\alpha_n, k_n} [f_\ell]   (x)  \\
     =&  \sum_{n=0 }^3  \psi_n     \frac{\partial }{\partial x_n} \left[  \frac{\sigma_n}{     
   \dfrac{d }{dx_n} e ( x_n^{\alpha_n})_{k _n}    }  I^{\beta_n}[f] (x) \right]  -   \sum_{n, \ell = 0 }^3  \psi_n \psi_\ell 
    T^{\sigma_n, \beta_n }_{\alpha_n, k_n} [f_\ell] (x).
\end{align*}
For each ${n=0,1,2,3}$ we see that    
\[\frac{\sigma_n}{\dfrac{d }{dx_n} e ( x_n^{\alpha_n})_{k _n}} I^{\beta_n}[f] (x) = \sum_{\ell=0}^3 \frac{\sigma_\ell}{\dfrac{d }{dx_\ell} e( x_\ell^{\alpha_\ell})_{k _\ell}}I^{\beta_\ell} [f](x) - \sum_{
   {
 \substack{\ell, n = 0 \\ \ell \ne n}}}^3 \frac{\sigma_\ell}{\dfrac{d }{dx_\ell} e( x_\ell^{\alpha_\ell})_{k _\ell}} I^{\beta_\ell} [f](x).\]
  
Therefore,
\begin{align*} 
& {}^{\psi}\mathcal D \left[ \sum_{\ell=0}^3  \frac{\sigma_\ell}{\dfrac{d }{dx_\ell} e( x_\ell^{\alpha_\ell})_{k _\ell}} I^{\beta_\ell} [f] \right](x)  \\ 
=  &  {}^{\psi}{\mathcal D}^{\sigma,\beta}_{\alpha, k} [f](x) + \sum_{n, \ell = 0 }^3 \psi_n \psi_\ell T^{\sigma_n, \beta_n }_{\alpha_n, k_n} [f_\ell] (x)  + 
     \sum_{
   {
      \substack{\ell, n = 0 \\ \ell \ne n}}}^3 \frac{\sigma_\ell}{     
\dfrac{d }{dx_\ell} e ( x_\ell^{\alpha_\ell})_{k _\ell}} \psi_n \frac{\partial }{\partial x_n} I^{\beta_\ell} [f](x).
\qedhere
\end{align*}
\end{proof}

Using the notation in Definition \ref{def100} and  Proposition \ref{prop1234} we have the following two corollaries:
\begin{corollary}\label{partcases1}Suppose that  
\[
\sigma = (\sigma_0,\sigma_1,\sigma_2,\sigma_3) \in (0,1]^4
.\]
 If ${\beta =(1,1,1,1)}$, then  
   \begin{align*} 
   {}^{\psi}{\mathcal D}^{\sigma,\beta}_{\alpha, k} [f](x)    =& \sum_{n, \ell = 0 }^3  \psi_n \psi_\ell \left(
  (1-\sigma_n)f_{\ell} (x) +   \sigma_n \frac{ \dfrac{\partial   f_\ell     }{\partial x_n}  (x)  }{  \dfrac{d }{dx_n} e ( x_n^{\alpha_n})_{k _n}    }\right) ,\\
  h_{n  }(x) = H^{\sigma_n, 1 }_{\alpha_n, k_n}[f_\ell](x) = & \frac{1-\sigma_{n}}{\sigma_n} \left[e(x_n^{\alpha_n})_{k _n} -1 \right],\\
  T^{\sigma_n, 1 }_{\alpha_n, k_n} [f_\ell] (x)= & 
      \frac{\partial }{\partial x_n} \left(  \frac{\sigma_n}{e^{h_{n}(x)   }    
   \dfrac{d }{dx_n} e ( x_n^{\alpha_n})_{k _n}    }  \right)     e^{h_{n}(x)}   f_\ell  (x)    ,\\
  I^{1}[f]   = & f, \\
          {}^{\psi}W_{\alpha, k}^{\sigma, \beta} [ f ] (x)   
 =  &
     \sum_{
   {
      \begin{array}{c}
        \ell =0=n \\
     \ell \neq n  
      \end{array}} 
   }^3 \frac{\sigma_\ell}{\dfrac{d }{dx_\ell} e ( x_\ell^{\alpha_\ell})_{k _\ell}} \psi_n \frac{\partial }{\partial x_n} f(x),  
\end{align*}
for all ${x \in \Omega}$ and \eqref{equa10} becomes 
\begin{align*}
{}^{\psi}\mathcal D \left[ \sum_{\ell=0}^3  \frac{\sigma_\ell}{     
   \dfrac{d }{dx_\ell} e ( x_\ell^{\alpha_\ell})_{k _\ell}    }  f    \right] (x)  
   =   & 
   {}^{\psi}{\mathcal D}^{\sigma,\beta}_{\alpha, k} [f](x)  +  A(x)  f(x) + {}^{\psi}W_{\alpha, k}^{\sigma, \beta} [ f ] (x) ,   
	\end{align*}
where 
\[A(x)= \sum_{n=0  }^3 \psi_n \frac{\partial }{\partial x_n} \left(\frac{\sigma_n}{e^{h_{n}(x)} \dfrac{d }{dx_n} e(x_n^{\alpha_n})_{k _n}} \right) e^{h_{n}(x)},\]
for all ${x\in \Omega}$.

Other important cases are as follows:
\begin{enumerate}
\item If ${\beta =(1,1,1,1)}$ and ${k=(1,1,1,1)}$, then  
\begin{align*} 
{}^{\psi}{\mathcal D}^{\sigma,\beta}_{\alpha, k} [f](x)    =& \sum_{n, \ell = 0 }^3  \psi_n \psi_\ell \left(
  (1-\sigma_n)f_{\ell} (x) +   \sigma_n \frac{ \dfrac{\partial   f_\ell     }{\partial x_n}  (x)  }{  {\alpha_n}  x_n^{\alpha_n-1}    }\right) ,\\
h_{n  }(x) = H^{\sigma_n, 1 }_{\alpha_n, 1}[f_\ell](x )   =  &  \frac{1-\sigma_{n}}{\sigma_n}     x_n^{\alpha_n}   ,\\
  T^{\sigma_n, 1 }_{\alpha_n, 1} [f_\ell] (x)= & 
      \frac{\partial }{\partial x_n} \left(  \frac{\sigma_n}{e^{h_{n}(x)   }    
   \alpha_n x_n^{\alpha_n-1}    }  \right)     e^{h_{n}(x)}   f_\ell  (x)    ,\\
  I^{1}[f]   = & f, \\
          {}^{\psi}W_{\alpha, k}^{\sigma, \beta} [ f ] (x)   
 =  &
     \sum_{
   {
      \substack{\ell, n = 0 \\ \ell \ne n}} 
   }^3 \frac{\sigma_\ell}{\alpha_\ell  x_\ell^{\alpha_\ell-1}}\psi_n \frac{\partial }{\partial x_n} f(x),  
 \end{align*}
 for all ${x \in \Omega}$ and \eqref{equa10} becomes   
\begin{align*}
{}^{\psi}\mathcal D \left[ \sum_{\ell=0}^3  \frac{\sigma_\ell}{{\alpha_\ell} x_\ell^{\alpha_\ell-1}} f \right] (x)  
   =   & {}^{\psi}{\mathcal D}^{\sigma,\beta}_{\alpha, k} [f](x)  +  A(x) f (x) + {}^{\psi}W_{\alpha, k}^{\sigma, \beta} [f] (x),
\end{align*}
where 
\[A(x)= \sum_{n=0}^3 \psi_n \frac{\partial}{\partial x_n} \left(\frac{\sigma_n}{e^{h_{n}(x)}\alpha_n x_n^{\alpha_n-1}}\right) e^{h_{n}(x)},\]
for all ${x\in \Omega}$.
\item If ${\beta =(1,1,1,1)}$ and ${k=(\infty, \infty, \infty ,\infty)}$, then  
   \begin{align*} 
   {}^{\psi}{\mathcal D}^{\sigma,\beta}_{\alpha, k} [f](x)    =& \sum_{n, \ell = 0 }^3  \psi_n \psi_\ell \left(
  (1-\sigma_n)f_{\ell} (x) +   \sigma_n \frac{ \dfrac{\partial   f_\ell     }{\partial x_n}  (x)  }{ \alpha_n x_n^{\alpha_n-1} e ^{x_n^{\alpha_n}}   }\right) ,\\
  h_{n  }(x) = H^{\sigma_n, 1 }_{\alpha_n, \infty }[f_\ell](x )   =  &  \frac{1-\sigma_{n}}{\sigma_n}  \left[   e^{ x_n^{\alpha_n}} -1 \right],\\
  T^{\sigma_n, 1 }_{\alpha_n, \infty } [f_\ell] (x)= & 
      \frac{\partial }{\partial x_n} \left(  \frac{\sigma_n}{e^{h_{n}(x)   }    
   \alpha_n x_n^{\alpha_n-1} e ^{x_n^{\alpha_n}}    }  \right)     e^{h_{n}(x)}   f_\ell  (x)    ,\\
  I^{1}[f]   = & f,  \\
          {}^{\psi}W_{\alpha, k}^{\sigma, \beta} [ f ] (x)   
 =  &
     \sum_{
   {
\substack{\ell, n = 0 \\ \ell \ne n}} 
   }^3 \frac{\sigma_\ell}{     
\alpha_\ell x_\ell^{\alpha_\ell-1} e^{x_\ell^{\alpha_\ell}}} \psi_n \frac{\partial }{\partial x_n} f(x),  
\end{align*}
for all ${x \in \Omega}$ and \eqref{equa10} becomes 
\begin{align*}
{}^{\psi}\mathcal D \left[ \sum_{\ell=0}^3  \frac{\sigma_\ell}{     
  \alpha_\ell x_\ell^{\alpha_\ell-1} e^ {x_\ell^{\alpha_\ell}} }f \right] (x)  
   =   & {}^{\psi}{\mathcal D}^{\sigma,\beta}_{\alpha, k} [f](x)  +  A(x)  f   (x)     + {}^{\psi}W_{\alpha, k}^{\sigma, \beta} [ f ] (x) ,   
\end{align*}
where 
\[ A(x)= \sum_{n=0  }^3  \psi_n   
    \frac{\partial }{\partial x_n} \left(  \frac{\sigma_n}{e^{h_{n}(x)   }    
   \alpha_n x_n^{\alpha_n-1} e^ {x_n^{\alpha_n}}    }  \right)     e^{h_{n}(x)} ,\]
for all ${x\in \Omega}$.
\end{enumerate}
\end{corollary}

\begin{corollary}\label{partcases2}Suppose  that  
\[
\rho = (\rho_0,\rho_1,\rho_2,\rho_3)\in (0,1]^4
.\]
 If ${\delta =(1,1,1,1)}$, then  
\begin{align*} 
{}^{\psi}{\mathcal D}^{\rho,\delta}_{r, \gamma, m} [g](x)  = &  
\sum_{n, \ell = 0 }^3   \psi_\ell \frac{ {\partial }^{\rho_n,1} g_\ell}{\partial x _{\gamma_n, m_n} } (x) \psi_n, \\
j_{n }(x) = H^{\rho_n, 1}_{\gamma_n, k_n}[g_{\ell}](x) =  & \frac{1-\rho_{n}}{\rho_n} \left[e ( x_n^{\gamma_n})_{m _n}- 1 \right], \\ 
  S^{\rho_n, 1 }_{\gamma_n, m_n} [g_\ell] (x) = & 
      \frac{\partial }{\partial x_n} \left(  \frac{\rho_n} { e^{j_{n }(x)   }    
   \dfrac{d }{dx_n} e ( x_n^{\gamma_n})_{m _n}    }  \right)     e^{j_{n }(x)}   g_\ell  (x)   ,\\
         {}^{\psi}V_{\gamma, m}^{\rho, \delta} [g ] (x)   =  & 
       \sum_{
   {
      \begin{array}{c}
        \ell =0=n \\
     \ell \neq n  
      \end{array}} 
   }^3       \frac{\rho_\ell}{     
   \dfrac{d }{dx_\ell} e ( x_\ell^{\gamma_\ell})_{m _\ell}    }      \frac{\partial }{\partial x_n}    g(x)  \psi_n ,
     \end{align*}
  and  identity   \eqref{equa11} is 
\begin{align*} 
  {}^{\psi}\mathcal D_r \left[ \sum_{\ell=0}^3  \frac{\rho_\ell}{     
   \dfrac{d }{dx_\ell} e ( x_\ell^{\gamma_\ell})_{m _\ell}    }   g    \right] (x)    = 
      {}^{\psi}{\mathcal D}^{\rho,\delta}_{r, \gamma, m} [g](x)  +  g  (x)  B(x) + 
{}^{\psi}V_{\gamma, m}^{\rho, \delta} [g ] (x),   \end{align*}
where 
\[B(x)=  \sum_{n=0  }^3       \frac{\partial }{\partial x_n} \left(  \frac{\rho_n} { e^{j_{n }(x)   }    
   \dfrac{d }{dx_n} e ( x_n^{\gamma_n})_{m _n}    }  \right)     e^{j_{n }(x)}       \psi_n,\]
for all ${x\in \Omega}$.
\begin{enumerate}
\item If ${\delta =(1,1,1,1)}$, then   
\begin{align*} 
{}^{\psi}{\mathcal D}^{\rho,\delta}_{r, \gamma, m} [g](x)  = &  
\sum_{n, \ell = 0 }^3 \psi_\ell \frac{ {\partial }^{\rho_n,1} g_\ell  }{ \partial  x _{\gamma_n, m_n} }   (x) \psi_n, \\
j_{n }(x ) = H^{\rho_n, 1 }_{\gamma_n, k_n}[g_{\ell}](x) = & \frac{1-\rho_{n}}{\rho_n} \left[e ( x_n^{\gamma_n})_{m _n}- 1 \right], \\ 
S^{\rho_n, 1 }_{\gamma_n, m_n} [g_\ell] (x) = & \frac{\partial }{\partial x_n} \left(\frac{\rho_n} {e^{j_{n }(x)}    
\dfrac{d }{dx_n} e ( x_n^{\gamma_n})_{m _n}} \right) e^{j_{n }(x)} g_\ell  (x), \\
{}^{\psi}V_{\gamma, m}^{\rho, \delta} [g](x) =  & 
       \sum_{
   {
\substack{\ell, n = 0 \\ \ell \ne n}} 
   }^3       \frac{\rho_\ell}{     
   \dfrac{d }{dx_\ell} e ( x_\ell^{\gamma_\ell})_{m _\ell}    }      \frac{\partial }{\partial x_n}    g(x)  \psi_n ,
\end{align*}
and  identity   \eqref{equa11} is 
\begin{align*} 
  {}^{\psi}\mathcal D_r \left[ \sum_{\ell=0}^3  \frac{\rho_\ell}{     
   \dfrac{d }{dx_\ell} e ( x_\ell^{\gamma_\ell})_{m _\ell}    }   g    \right] (x)    = 
      {}^{\psi}{\mathcal D}^{\rho,\delta}_{r, \gamma, m} [g](x)  +  g  (x)  B(x) + 
{}^{\psi}V_{\gamma, m}^{\rho, \delta} [g ] (x),   \end{align*}
where 
\[B(x)=  \sum_{n=0  }^3       \frac{\partial }{\partial x_n} \left(  \frac{\rho_n} { e^{j_{n }(x)   }    
   \dfrac{d }{dx_n} e ( x_n^{\gamma_n})_{m _n}    }  \right)     e^{j_{n }(x)}       \psi_n,\]
for all ${x\in \Omega}$.
\item If ${\delta =(1,1,1,1)}$ and ${m= (1,1,1,1)}$, then  
\begin{align*} 
{}^{\psi} {\mathcal D}^{\rho,\delta}_{r, \gamma, m} [g](x) = & \sum_{n, \ell = 0 }^3 \psi_\ell \frac{ {\partial }^{\rho_n,1} g_\ell  }{ \partial  x _{\gamma_n, 1} } (x) \psi_n, \\
j_{n }(x ) =   H^{\rho_n, 1 }_{\gamma_n, k_n}[g_{\ell}](x )   =  & \frac{1-\rho_{n}}{\rho_n}     x_n^{\gamma_n}   , \\ 
  S^{\rho_n, 1 }_{\gamma_n, 1} [g_\ell] (x) = & 
      \frac{\partial }{\partial x_n} \left(  \frac{\rho_n} { e^{j_{n }(x)   }    
    \gamma_n  x_n^{\gamma_n-1}    }  \right)     e^{j_{n }(x)}   g_\ell  (x)   ,\\
         {}^{\psi}V_{\gamma, m}^{\rho, \delta} [g ] (x)   =  & 
       \sum_{
   {
      \substack{\ell, n = 0 \\ \ell \ne n}} 
   }^3       \frac{\rho_\ell}{     
   \gamma_\ell  x_\ell^{\gamma_\ell-1}} \frac{\partial }{\partial x_n} g(x)  \psi_n ,
\end{align*}
and  identity   \eqref{equa11} is 
\begin{align*} 
  {}^{\psi}\mathcal D_r \left[ \sum_{\ell=0}^3  \frac{\rho_\ell}{     
  \gamma_\ell x_\ell^{\gamma_\ell-1}      }   g    \right] (x)    = 
      {}^{\psi}{\mathcal D}^{\rho,\delta}_{r, \gamma, m} [g](x)  +  g  (x)  B(x) + 
{}^{\psi}V_{\gamma, m}^{\rho, \delta} [g ] (x),   \end{align*}
where 
\[B(x) = \sum_{n=0}^3 \frac{\partial}{\partial x_n} \left(\frac{\rho_n} {e^{j_{n }(x)} \gamma_n x_n^{\gamma_n-1}}\right) e^{j_{n }(x)}  \psi_n,\]
for all ${x\in \Omega}$.
\item If ${\delta =(1,1,1,1)}$ and ${m= (\infty , \infty,\infty,\infty)}$, then   
\begin{align*} 
{}^{\psi}{\mathcal D}^{\rho,\delta}_{r, \gamma, m} [g](x)  = &  
\sum_{n, \ell = 0 }^3   \psi_\ell \frac{ {\partial }^{\rho_n,1} g_\ell  }{ \partial  x _{\gamma_n, \infty}}(x) \psi_n, \\
j_{n }(x ) =   H^{\rho_n, 1 }_{\gamma_n, k_n}[g_{\ell}](x )   =  & \frac{1-\rho_{n}}{\rho_n} \left[  e^{x_n^{\gamma_n}}   -1 \right] , \\ 
  S^{\rho_n, 1 }_{\gamma_n, \infty } [g_\ell] (x) = & 
      \frac{\partial }{\partial x_n} \left(  \frac{\rho_n} { e^{j_{n }(x)   }    
   \gamma_n x_n^{\gamma_n-1} e^{x_n^{\gamma_n}}    }  \right)     e^{j_{n }(x)}   g_\ell  (x)   ,\\
         {}^{\psi}V_{\gamma, m}^{\rho, \delta} [g ] (x)   =  & 
       \sum_{
   {
     \substack{\ell, n = 0 \\ \ell \ne n}} 
   }^3       \frac{\rho_\ell}{     
   \gamma_\ell x_\ell^{\gamma_\ell-1} e^{x_\ell^{\gamma_\ell}}  
      }      \frac{\partial }{\partial x_n}    g(x)  \psi_n ,
\end{align*}
and  identity   \eqref{equa11} is 
\begin{align*}
  {}^{\psi}\mathcal D_r \left[ \sum_{\ell=0}^3  \frac{\rho_\ell}{     
     \gamma_\ell x_\ell^{\gamma_\ell-1} e^{x_\ell^{\gamma_\ell}} 
   }   g    \right] (x)    = 
      {}^{\psi}{\mathcal D}^{\rho,\delta}_{r, \gamma, m} [g](x)  +  g  (x)  B(x) + 
{}^{\psi}V_{\gamma, m}^{\rho, \delta} [g ] (x),   \end{align*}
where 
\[B(x)=  \sum_{n=0  }^3       \frac{\partial }{\partial x_n} \left(  \frac{\rho_n} { e^{j_{n }(x)   }    
     \gamma_n x_n^{\gamma_n-1} e^{x_n^{\gamma_n}} 
       }  \right)     e^{j_{n }(x)}       \psi_n,\]
for all ${x\in \Omega}$.
\end{enumerate}
\end{corollary}
 
Note that if ${\beta=(1,1,1,1)}$, then the condition ${f(\Omega) \subset \mathbb H^{+}}$
is not necessary.
 
The following theorem presents Borel--Pompeiu and Stokes formulas induced by the operators
${}^{\psi}{\mathcal D}^{\sigma,\beta}_{\alpha, k}$ and ${}^{\psi}{\mathcal D}^{\rho,\delta}_{r, \gamma, m}$.

\begin{theorem}\label{Borel-Pompeiu and Stokes formulas}
Let ${\Omega\subset \mathbb H^+}$ be a domain such that ${\partial \Omega\subset \mathbb H^+}$ is a 3-dimensional compact smooth surface. Using the notation in Definition \ref{def100},   let 
\[
f= \sum_{\ell=0}^3 \psi_\ell  f_\ell , \quad  g= \sum_{\ell=0}^3 \psi_\ell  g_\ell \in C^1(\Omega,\mathbb H)
,\]
 where ${f_\ell, g_\ell}$ are ${\mathbb{R}^{+}}$-valued functions. Suppose that  
\[
\sigma = (\sigma_0,\sigma_1,\sigma_2,\sigma_3), \ \rho = (\rho_0,\rho_1,\rho_2,\rho_3)\in (0,1]^4
.\]

Then
\begin{align}\label{BorelHypGral}  &  \int_{\partial \Omega} K_{\psi}(\tau-x)\sigma_{\tau}^{\psi} \left(\sum_{\ell=0}^3  \frac{\sigma_\ell}{\dfrac{d }{d\tau_\ell} e ( \tau_\ell^{\alpha_\ell})_{k _\ell}}I^{\beta_\ell} [f] (\tau) \right) \nonumber \\
& + \int_{\partial \Omega}  \left(\sum_{\ell=0}^3  \frac{\rho_\ell}{\dfrac{d }{d\tau_\ell} e ( \tau_\ell^{\gamma_\ell})_{m _\ell}} I^{\delta_\ell} [g] (\tau) \right) \sigma_{\tau}^{\psi} K_{\psi}(\tau-x)  \nonumber  \\ 
& - \int_{\Omega} \left[  K_{\psi} (y-x)\,  {}^{\psi}{\mathcal D}^{\sigma,\beta}_{\alpha, k} [f](y)   + 
     {}^{\psi}{\mathcal D}^{\rho,\delta}_{r, \gamma, m} [g](y) \,K_{\psi} (y-x) \right]dy   \nonumber \\
& - \int_{\Omega}  K_{\psi} (y-x)  \left[   \sum_{n, \ell = 0 }^3  \psi_n \psi_\ell  T^{\sigma_n, \beta_n }_{\alpha_n, k_n} [f_\ell] (y)  + 
{}^{\psi}W_{\alpha, k}^{\sigma, \beta} [ f ] (y) \right] dy \nonumber \\
& - \int_{\Omega} \left[\sum_{n, \ell = 0 }^3 \psi_\ell S^{\rho_n, \delta_n }_{\gamma_n, m_n} [g_\ell] (y)  \psi_n + 
{}^{\psi}V_{\gamma, m}^{\rho, \delta} [g ] (y) \right] K_{\psi} (y-x) dy  \nonumber \\
=  &  \left\{ \begin{array}{ll} \begin{aligned} &  \displaystyle  \sum_{\ell=0}^3  \frac{\sigma_\ell}{\dfrac{d }{dx_\ell} e ( x_\ell^{\alpha_\ell})_{k _\ell}} I^{\beta_\ell} [f] (x) + \sum_{\ell=0}^3 \frac{\rho_\ell}{\dfrac{d}{dx_\ell} e (x_\ell^{\gamma_\ell})_{m _\ell}} I^{\delta_\ell} [g]  
(x) \end{aligned}, &  x\in \Omega,  \\ 0 , &  x\in \mathbb H\setminus\overline{\Omega}.\end{array} \right. 
\end{align} 
In addition, 
\begin{align}\label{stokefrac}
    \int_{\partial \Omega} & \left(  \sum_{\ell=0}^3  \frac{\rho_\ell}{     
   \dfrac{d }{dx_\ell} e ( x_\ell^{\gamma_\ell})_{m _\ell}    }    I^{\delta_\ell} [g]  
  \right)    \sigma^\psi_x \left(  \sum_{\ell=0}^3  \frac{\sigma_\ell}{     
   \dfrac{d }{dx_\ell} e ( x_\ell^{\alpha_\ell})_{k _\ell}    }    I^{\beta_\ell} [f]  (x)  
  \right)              \nonumber      \\
	=  &   \int_{\Omega } \left( g\, {}^{\psi}{\mathcal D}^{\sigma,\beta}_{\alpha, k} [f](x) 
	+  {}^{\psi}{\mathcal D}^{\rho,\delta}_{r, \gamma, m} [g](x) \,f  \right) dx \nonumber \\
 &	   +  \int_{\Omega }  g(x) \left[  \sum_{n, \ell = 0 }^3  \psi_n \psi_\ell 
    T^{\sigma_n, \beta_n }_{\alpha_n, k_n} [f_\ell] (x)  + 
{}^{\psi}W_{\alpha, k}^{\sigma, \beta} [ f ] (x)\right] dx  \nonumber  \\
 & +
	\int_{\Omega}
 \left[    \sum_{n, \ell = 0 }^3   \psi_\ell 
    S^{\rho_n, \delta_n }_{\gamma_n, m_n} [g_\ell] (x)  \psi_n + 
{}^{\psi}V_{\gamma, m}^{\rho, \delta} [g ] (x) \right]  f(x) dx.
\end{align}
\end{theorem}

\begin{proof}
This follows directly from Definition \ref{def100}  using functions

\[
\displaystyle   \sum_{\ell=0}^3  \frac{\sigma_\ell}{\dfrac{d }{dx_\ell} e ( x_\ell^{\alpha_\ell})_{k _\ell}} I^{\beta_\ell} [f](x),
\]
 and 
\[
\displaystyle \sum_{\ell=0}^3  \frac{\rho_\ell}{\dfrac{d }{dx_\ell} e ( x_\ell^{\gamma_\ell})_{m _\ell}} I^{\delta_\ell} [g](x),
\]
 and identities  \eqref{equa10} and \eqref{equa11} in formulas \eqref{BorelHyp} and \eqref{StokesHyp}.
\end{proof}

\begin{remark}
In formulas \eqref{BorelHypGral} and \eqref{stokefrac},  the operators  
\begin{center}
${}^{\psi}{\mathcal D}^{\sigma,\beta}_{\alpha, k}$
 and  
${}^{\psi}{\mathcal D}^{\rho,\delta}_{r, \gamma, m}$
\end{center}
 reflect the phenomenon of duality in quaternionic  analysis due to the non-commutativity of quaternionic algebra.
\end{remark}

\begin{corollary}[Cauchy-type Formula] \label{CauchyFRactal} 
Let ${\Omega\subset \mathbb H^+}$ be a domain such that ${\partial \Omega\subset \mathbb H^+}$ is a 3-dimensional compact smooth surface. Suppose that  
\[
\sigma = (\sigma_0,\sigma_1,\sigma_2,\sigma_3), \ \rho = (\rho_0,\rho_1,\rho_2,\rho_3)\in (0,1]^4
.\]
 If 
 
\[
f\in        C^1(\Omega, \mathbb H)  \cap \textrm{Ker} {}^{\psi}{\mathcal D}^{\sigma,\beta}_{\alpha, k},
\]
 and 
\[
g \in 
C^1  (\Omega, \mathbb H) \cap \textrm{Ker} {}^{\psi}{\mathcal D}^{\rho,\delta}_{r, \gamma, m}
,\]
  then 
\begin{align*}  
\int_{\partial \Omega} & K_{\psi}(\tau-x)\sigma_{\tau}^{\psi} \left(\sum_{\ell=0}^3  \frac{\sigma_\ell}{     
   \dfrac{d }{d\tau_\ell} e ( \tau_\ell^{\alpha_\ell})_{k _\ell}} I^{\beta_\ell} [f](\tau) \right) \nonumber \\
    &  + \int_{\partial \Omega}  \left(\sum_{\ell=0}^3  \frac{\rho_\ell}{\dfrac{d }{d\tau_\ell} e (\tau_\ell^{\gamma_\ell})_{m _\ell}} I^{\delta_\ell} [g](\tau)\right) \sigma_{\tau}^{\psi} K_{\psi}(\tau-x)\nonumber \\
& - \int_{\Omega}  K_{\psi} (y-x)  \left[   \sum_{n, \ell = 0 }^3  \psi_n \psi_\ell T^{\sigma_n, \beta_n }_{\alpha_n, k_n} [f_\ell] (y)  + 
{}^{\psi}W_{\alpha, k}^{\sigma, \beta} [ f ] (y)\right] dy \nonumber \\
& - \int_{\Omega} \left[\sum_{n, \ell = 0 }^3   \psi_\ell S^{\rho_n, \delta_n }_{\gamma_n, m_n} [g_\ell] (y)  \psi_n + 
{}^{\psi}V_{\gamma, m}^{\rho, \delta} [g](y)\right] K_{\psi} (y-x) dy   \nonumber \\
		=  &  \left\{ \begin{array}{ll} \begin{aligned} &  \displaystyle  \sum_{\ell=0}^3  \frac{\sigma_\ell}{     
   \dfrac{d }{dx_\ell} e ( x_\ell^{\alpha_\ell})_{k _\ell}} I^{\beta_\ell} [f](x) +  \sum_{\ell=0}^3  \frac{\rho_\ell}{     
   \dfrac{d }{dx_\ell} e ( x_\ell^{\gamma_\ell})_{m _\ell}} I^{\delta_\ell} [g](x) \end{aligned}, &  x\in \Omega,  \\ 0 , &  x\in \mathbb H\setminus\overline{\Omega},                     
\end{array} \right. 
\end{align*} 
and 
\begin{align*}
\int_{\partial \Omega}&\left(\sum_{\ell=0}^3 \frac{\rho_\ell}{\dfrac{d }{dx_\ell} e(x_\ell^{\gamma_\ell})_{m _\ell}} I^{\delta_\ell} [g]\right) \sigma^\psi_x \left(  \sum_{\ell=0}^3  \frac{\sigma_\ell}{\dfrac{d }{dx_\ell} e (x_\ell^{\alpha_\ell})_{k _\ell}} I^{\beta_\ell} [f](x) \right) \nonumber \\
=  &  \int_{\Omega }  g(x) \left[  \sum_{n, \ell = 0 }^3  \psi_n \psi_\ell T^{\sigma_n, \beta_n }_{\alpha_n, k_n} [f_\ell] (x)  + 
{}^{\psi}W_{\alpha, k}^{\sigma, \beta} [ f ] (x)\right] dx  \nonumber  \\
 & + \int_{\Omega} \left[\sum_{n, \ell = 0 }^3 \psi_\ell S^{\rho_n, \delta_n }_{\gamma_n, m_n} [g_\ell] (x)  \psi_n + 
{}^{\psi}V_{\gamma, m}^{\rho, \delta} [g ] (x) \right]f(x)dx.
\end{align*}
\end{corollary}

\begin{corollary}\label{CauchyFRactal1}
Suppose that ${\beta=(1,1,1,1)}$, ${\delta = (1,1,1,1)}$ and   
\[
\sigma = (\sigma_0,\sigma_1,\sigma_2,\sigma_3), \ \rho = (\rho_0,\rho_1,\rho_2,\rho_3)\in (0,1]^4
.\]
Let ${\Omega\subset \mathbb H^+}$ be a domain such that ${\partial \Omega\subset \mathbb H^+}$ is a 3-dimensional compact smooth surface. Let 
\[
f= \sum_{\ell=0}^3 \psi_\ell  f_\ell \in  C^1(\Omega, \mathbb H)  \cap \textrm{Ker} {}^{\psi}{\mathcal D}^{\sigma,\beta}_{\alpha, k}
,\]

\[
g= \sum_{\ell=0}^3 \psi_\ell  g_\ell \in C^1  (\Omega, \mathbb H) \cap \textrm{Ker} {}^{\psi}{\mathcal D}^{\rho,\delta}_{r, \gamma, m}
,\]
 where ${f_\ell, g_\ell}$ are ${\mathbb{R}}$-valued functions. 
\begin{enumerate}
\item If ${k=(1,1,1,1)}$ and ${m= (1,1,1,1)}$,  then 
\begin{align*}  \int_{\partial \Omega} & K_{\psi}(\tau-x)\sigma_{\tau}^{\psi} \left(\sum_{\ell=0}^3  \frac{\sigma_\ell}{     
    \alpha_\ell \tau_\ell^{\alpha_\ell- 1 }    }     \right)  f   (\tau)   \\ &+ \int_{\partial \Omega}  g (\tau)  \left(
    \sum_{\ell=0}^3  \frac{\rho_\ell}{     
   \gamma_\ell \tau_\ell^{\gamma_\ell-1}    }   
    \right)   \sigma_{\tau}^{\psi} K_{\psi}(\tau-x)  \nonumber  \\ 
&
- \int_{\Omega}  K_{\psi} (y-x)  \left[   \sum_{n, \ell = 0 }^3  \psi_n \psi_\ell 
    T^{\sigma_n, 1 }_{\alpha_n, 1} [f_\ell] (y)  + 
{}^{\psi}W_{\alpha, k}^{\sigma, \beta} [ f ] (y)   \right]  dy
 \nonumber \\
& 
- \int_{\Omega}   
 \left[  
    \sum_{n, \ell = 0 }^3   \psi_\ell 
    S^{\rho_n, 1 }_{\gamma_n, 1} [g_\ell] (y)  \psi_n + 
{}^{\psi}V_{\gamma, m}^{\rho, \delta} [g ] (y)
  \right] K_{\psi} (y-x) dy   \nonumber \\
		=  &  \left\{ \begin{array}{ll}   \displaystyle   f (x) \displaystyle   \sum_{\ell=0}^3  \frac{\sigma_\ell}{     
  \alpha_\ell  x_\ell^{\alpha_\ell-1}} + g(x) \sum_{\ell=0}^3  \frac{\rho_\ell}{\gamma_\ell x_\ell^{\gamma_\ell-1}}, &  x\in \Omega, \\ 0 , &  x\in \mathbb H\setminus\overline{\Omega},                     
\end{array} \right. 
\end{align*} 
and 
\begin{align*}
	\int_{\partial \Omega} &  g  (x) \left(  \sum_{\ell=0}^3  \frac{\rho_\ell}{     
    \gamma_\ell x_\ell^{\gamma_\ell-1}    }  
  \right)    \sigma^\psi_x \left(  \sum_{\ell=0}^3  \frac{\sigma_\ell}{     
   \alpha_\ell x_\ell^{\alpha_\ell-1}    }   \right)   f  (x)  
              \nonumber      \\
	=  &  \int_{\Omega }  g(x) \left[  \sum_{n, \ell = 0 }^3  \psi_n \psi_\ell 
    T^{\sigma_n, 1 }_{\alpha_n, 1} [f_\ell] (x)  + 
{}^{\psi}W_{\alpha, k}^{\sigma, \beta} [ f ] (x)\right] dx  \nonumber  \\
 & +
	\int_{\Omega}
 \left[    \sum_{n, \ell = 0 }^3   \psi_\ell 
    S^{\rho_n, 1 }_{\gamma_n, 1} [g_\ell] (x)  \psi_n + 
{}^{\psi}V_{\gamma, m}^{\rho, \delta} [g ] (x) \right]  f(x) dx.	
\end{align*}
\item If ${k=(\infty, \infty , \infty ,\infty)}$ and ${m= (\infty, \infty , \infty ,\infty)}$,  then 
 \begin{align*}   \int_{\partial \Omega} & K_{\psi}(\tau-x)\sigma_{\tau}^{\psi} \left(\sum_{\ell=0}^3  \frac{\sigma_\ell}{     
  \alpha_\ell  \tau_\ell^{\alpha_\ell-1}  e^{\tau_\ell^{\alpha_\ell}}    }     \right) f (\tau)  \\ &+ \int_{\partial \Omega} g(\tau) \left(
    \sum_{\ell=0}^3  \frac{\rho_\ell}{     
      \gamma_\ell  \tau_\ell^{\gamma_\ell-1}  e^{\tau_\ell^{\gamma_\ell}}
      }     
    \right)   \sigma_{\tau}^{\psi} K_{\psi}(\tau-x)  \nonumber  \\ 
&
- \int_{\Omega}  K_{\psi} (y-x)  \left[   \sum_{n, \ell = 0 }^3  \psi_n \psi_\ell 
    T^{\sigma_n, 1 }_{\alpha_n, \infty} [f_\ell] (y)  + 
{}^{\psi}W_{\alpha, k}^{\sigma, \beta} [ f ] (y)   \right]  dy
 \nonumber \\
& 
- \int_{\Omega}   
 \left[  
    \sum_{n, \ell = 0 }^3   \psi_\ell 
    S^{\rho_n,  1 }_{\gamma_n, \infty} [g_\ell] (y)  \psi_n + 
{}^{\psi}V_{\gamma, m}^{\rho, \delta} [g ] (y)
  \right] 
 K_{\psi} (y-x) dy   \nonumber \\
		=  &  \left\{ \begin{array}{ll}  \displaystyle  f(x)  \displaystyle  \sum_{\ell=0}^3  \frac{\sigma_\ell}{ 
		\alpha_\ell  x_\ell^{\alpha_\ell-1}  e^{x_\ell^{\alpha_\ell}}
   }    + g(x) \sum_{\ell=0}^3  \frac{\rho_\ell}{  
      \gamma_\ell  x_\ell^{\gamma_\ell-1}  e^{x_\ell^{\gamma_\ell}}   
     }, &  x\in \Omega,  \\ 0 , &  x\in \mathbb H\setminus\overline{\Omega},                    
\end{array} \right. 
\end{align*} 
and 
\begin{align*}
	\int_{\partial \Omega} & g(x) \left(  \sum_{\ell=0}^3  \frac{\rho_\ell}{     
  	\gamma_\ell  x_\ell^{\gamma_\ell-1}  e^{x_\ell^{\gamma_\ell}}     }
  \right)    \sigma^\psi_x \left(  \sum_{\ell=0}^3  \frac{\sigma_\ell}{     
  	\alpha_\ell  x_\ell^{\alpha_\ell-1}  e^{x_\ell^{\alpha_\ell}}
   }    
  \right)    f (x)          \nonumber      \\
	=  &  \int_{\Omega }  g(x) \left[  \sum_{n, \ell = 0 }^3  \psi_n \psi_\ell 
    T^{\sigma_n, 1 }_{\alpha_n, \infty } [f_\ell] (x)  + 
{}^{\psi}W_{\alpha, k}^{\sigma, \beta} [ f ] (x)\right] dx  \nonumber  \\
 & +
	\int_{\Omega}
 \left[    \sum_{n, \ell = 0 }^3   \psi_\ell 
    S^{\rho_n, 1 }_{\gamma_n, \infty } [g_\ell] (x)  \psi_n + 
{}^{\psi}V_{\gamma, m}^{\rho, \delta} [g ] (x) \right]  f(x) dx.
	\end{align*}
	\end{enumerate} 
	\end{corollary}

\begin{corollary}\label{CauchyFRactal2}
Suppose that
 ${\beta=(1,1,1,1)}$, ${\delta = (1,1,1,1)}$ and    
\[
\sigma = (\sigma_0,\sigma_1,\sigma_2,\sigma_3), \ \rho = (\rho_0,\rho_1,\rho_2,\rho_3)\in (0,1]^4
.\]
Let ${\Omega\subset \mathbb H^+}$ be a domain such that ${\partial \Omega\subset \mathbb H^+}$ is a 3-dimensional compact smooth surface. Let
\[
f= \sum_{\ell=0}^3 \psi_\ell  f_\ell \in  C^1(\Omega, \mathbb H)  \cap \textrm{Ker} {}^{\psi}{\mathcal D}^{\sigma,\beta}_{\alpha, k}
,\]

\[
g= \sum_{\ell=0}^3 \psi_\ell  g_\ell \in C^1  (\Omega, \mathbb H) \cap \textrm{Ker} {}^{\psi}{\mathcal D}^{\rho,\delta}_{r, \gamma, m}
,\]
 where ${f_\ell, g_\ell}$ are ${\mathbb{R}}$-valued functions.
\begin{enumerate}	
\item If ${k=(1,1,1,1)}$ and ${m=(\infty, \infty , \infty ,\infty)}$, then
  \begin{align*}   \int_{\partial \Omega} & K_{\psi}(\tau-x)\sigma_{\tau}^{\psi} \left(\sum_{\ell=0}^3  \frac{\sigma_\ell}{     
   \alpha_\ell \tau_\ell^{\alpha_\ell-1}   }   \right) f  (\tau)   \\ &+ \int_{\partial \Omega}   g(\tau) \left(
    \sum_{\ell=0}^3  \frac{\rho_\ell}{     
  \gamma_\ell \tau^{\gamma_\ell -1}  e ^{\tau_\ell^{\gamma_\ell} }    }    
    \right)   \sigma_{\tau}^{\psi} K_{\psi}(\tau-x)  \nonumber  \\ 
&
- \int_{\Omega}  K_{\psi} (y-x)  \left[   \sum_{n, \ell = 0 }^3  \psi_n \psi_\ell 
    T^{\sigma_n, 1 }_{\alpha_n, 1} [f_\ell] (y)  + 
{}^{\psi}W_{\alpha, k}^{\sigma, \beta} [ f ] (y)   \right]  dy
 \nonumber \\
& 
- \int_{\Omega}   
 \left[  
    \sum_{n, \ell = 0 }^3   \psi_\ell 
    S^{\rho_n, 1 }_{\gamma_n, \infty } [g_\ell] (y)  \psi_n + 
{}^{\psi}V_{\gamma, m}^{\rho, \delta} [g ] (y)
  \right] 
 K_{\psi} (y-x)  dy   \nonumber \\
		=  &  \left\{ \begin{array}{ll}   \displaystyle   f(x)\sum_{\ell=0}^3  \frac{\sigma_\ell}{     
  \alpha_\ell x_\ell^{\alpha_\ell-1}    }      + g(x)  \sum_{\ell=0}^3  \frac{\rho_\ell}{     
    \gamma_\ell x_\ell^{\gamma_\ell-1}  e ^{ x_\ell^{\gamma_\ell} }   }    
     , &  x\in \Omega,  \\ 0 , &  x\in \mathbb H\setminus\overline{\Omega},                     
\end{array} \right. 
\end{align*} 
and 
\begin{align*}
 \int_{\partial \Omega} & g(x) \left(  \sum_{\ell=0}^3  \frac{\rho_\ell}{     
  \gamma_\ell x_\ell^{\gamma_\ell-1}  e ^{ x_\ell^{\gamma_\ell} }     }     
  \right)    \sigma^\psi_x   \left(  \sum_{\ell=0}^3  \frac{\sigma_\ell}{     
   \alpha_\ell  x_\ell^{\alpha_\ell-1}    }      \right)   f(x)            \nonumber      \\
	=  &  \int_{\Omega }  g(x) \left[  \sum_{n, \ell = 0 }^3  \psi_n \psi_\ell 
    T^{\sigma_n, 1 }_{\alpha_n, 1} [f_\ell] (x)  + 
{}^{\psi}W_{\alpha, k}^{\sigma, \beta} [ f ] (x)\right] dx  \nonumber  \\
 & +
	\int_{\Omega}
 \left[    \sum_{n, \ell = 0 }^3   \psi_\ell 
    S^{\rho_n, 1 }_{\gamma_n, \infty } [g_\ell] (x)  \psi_n + 
{}^{\psi}V_{\gamma, m}^{\rho, \delta} [g ] (x) \right]  f(x) dx.
	\end{align*}
\item For ${k=(\infty, \infty , \infty ,\infty)}$ and ${m =(1,1,1,1)}$ an analogous result holds.
\end{enumerate} 
\end{corollary}

\section{Concluding remarks and future work}

This {work} establishes the foundations of a quaternionic function theory associated with  a proportional and {fractal} ${\psi}$-Fueter operator associated with a fractal measure. Also this {work} extends the quaternionic hyperholomorphic function  theory. So what other results can be extended to this recent function theory?

Since we are working in a quaternionic environment, some direct consequences and directions for future work include the following:
\begin{enumerate}
\item Proportional-fractal operators of divergence and rotation type, and proportional-fractal vector fields with respect to fractal measure functions, using convenient structural sets, among other important systems.

\item Proportional-fractal Cauchy--Riemann-type operators and  proportional-fractal holomorphic functions in two complex variables with respect  to fractal measure functions using convenient structural sets, among other important systems.

\item  The  Fueter operator satisfies  an important   property called  conformal covariance, closely related to conformal mappings, which we think could be extended to the  ${\beta}$-proportional fractal Fueter operator.

\item There is a wide range of quaternionic function  modules associated with the Fueter operator, such as hyperholomorphic Bergman spaces, Hardy spaces, Dirichlet spaces, among others, which could be generalized using   the quaternionic function module   induced by the   ${\beta}$-proportional fractal Fueter operator. 
\end{enumerate}

The computations presented in this article can be extended directly to the Clifford algebra framework to obtain a proportional-fractal Dirac-type operator and  proportional-fractal monogenic functions with respect  to fractal measure functions.  

The paper \cite{GMB} introduces a generalized fractal derivative with respect to a  truncated quaternionic exponential function on slices, to present numerous families of quaternionic Hilbert modules of slice regular functions, in which the Bergman and Dirichlet quaternionic modules constitute two important elements of one of these families.

Therefore, an interesting question arises: Does the proportional and fractal  ${\psi}$-Fueter operator, associated with a fractal measure, induce families of quaternionic Hilbert modules that extend the Bergman and Dirichlet quaternionic modules in quaternionic analysis?

\small{
    
}

\end{document}